\documentclass[12pt]{amsart}%
\usepackage{amsfonts}
\usepackage{amsmath}
\usepackage{soul}
\setstcolor{blue}
\sethlcolor{red}
\usepackage{color}
\usepackage{tikz,pgf}
\usetikzlibrary{arrows}
\usepackage{amssymb}
\usepackage{graphicx}%
\usepackage[all]{xy}\SelectTips{cm}{}
\definecolor{lightblue}{rgb}{.30,.95,.2}
\sethlcolor{lightblue}

\newcommand{\mb}[1]{\mathbf{ #1}}

\providecommand{\U}[1]{\protect\rule{.1in}{.1in}}
\newtheorem{theorem}{Theorem}
\theoremstyle{plain}

\newtheorem{corollary}{Corollary}

\newtheorem{definition}{Definition}
\newtheorem{example}{Example}

\newtheorem{lemma}{Lemma}

\newtheorem{proposition}{Proposition}
\newtheorem{remark}{Remark}

\numberwithin{equation}{section}
\begin{document}
\title{Representing quantum structures as near semirings}
\author[Bonzio S.]{S. Bonzio}
\address{Stefano Bonzio, University of Cagliari\\
Italy}
\email{stefano.bonzio@gmail.com}

\author[Chajda I.]{I. Chajda}
\address{Ivan Chajda, Palack\'{y} University, Olomouc\\
Czech Republic}
\email{ivan.chajda@upol.cz}

\author[Ledda A.]{A. Ledda}
\address{Antonio Ledda, University of Cagliari\\
Italy}
\email{antonio.ledda@unica.it}
%
%

\begin{abstract}

In this paper we introduce the notion of near semiring with involution. Generalizing the theory of semirings we aim at represent quantum structures, such as basic algebras and orthomodular lattices, in terms of near semirings with involution. In particular, after discussing several properties of near semirings, we introduce the so-called \L ukasiewicz near semirings, as a particular case of near semirings, and we show that every basic algebra is representable as (precisely, it is term equivalent to) a near semiring. In the particular case in which a \L ukasiewicz near semiring is also a semiring, we obtain as a corollary a representation of MV-algebras as semirings. Analogously, by introducing a particular subclass of \L ukasiewicz near semirings, that we termed orthomodular near semirings, we obtain a representation of orthomodular lattices. In the second part of the paper, we discuss several universal algebraic properties of \L ukasiewicz near semirings and we show that the variety of involutive integral near semirings is a Church variety. This yields a neat equational characterization of central elements of this variety. As a byproduct of such, we obtain several direct decomposition theorems for this class of algebras.
\end{abstract}
\thanks{{\bf Corresponding author}: {\bf Antonio Ledda}, {\ttfamily antonio.ledda@unica.it}}
\maketitle

\section{Introduction}
\noindent It is a long-dated result, due to Marshall Stone \cite{St36}, that the theory of Boolean algebras (the algebraic counterpart of classical logic) can be framed within the theory of rigs, through the concept of Boolean ring.
More recently, in the last decade, the relations between prominent algebraic structures from many-valued logics and (semi)ring theory have stirred a renewed attention (see, e.g., \cite{Gerla03, BdN09}).
It was shown by Belluce, Di Nola, Ferraioli \cite{DiNola13} and Gerla \cite{Gerla03} that MV-algebras (the algebraic semantics of infinite-valued \L ukasiewicz logic) can be viewed as particular semirings: MV-semirings. 
The results achieved in the MV-algebras' context are extremely interesting and promising. 
Quoting A. Di Nola and C. Russo \cite{DR15}:  
\begin{quotation}
[...] besides serving MV-algebra theory they suggest a possible Òpayback,Ó namely, that MV-algebras can on their turn give ideas and tools to semiring and semifield theories. It is worth noticing also that, as well as MV-algebras, various other logic-related algebraic structures can be viewed as special idempotent semirings, and therefore this approach could be further extended.
\end{quotation}

Taking up their suggestion, we will show in this paper that this method can be fruitfully raised to a considerably general level. Indeed, we will see that a number of algebraic structures of major importance to non classical logics are representable as semiring-like structures.

This paper will be mainly focused on basic algebras and orthomodular lattices. Basic algebras were introduced by R. Hala{\v{s}}, J. K{\"u}hr and one of the authors of this article as a common generalization of both MV-algebras --the algebraic alter-ego of \L ukasiewicz many-valued logic -- and orthomodular lattices (the interested reader may consult \cite{Chajda15} and \cite{Chajda07} for details)  -- the algebraic counterpart of the logic of quantum mechanics (for an extensive discussion we refer to \cite{Beran11, Kalmbach83}). Inspired by the results in \cite{DiNola13}, it seems natural to hunt for an appropriate notion that would play, in the wider domain of Basic algebras, the same role that MV-semirings interpret in the context of MV-algebras. We believe that this concept will shine some light on the theory both of MV-algebras and orthomodular lattices. It will be interesting, in fact, to examine how this general notion specifies to context so far apart from each other. Indeed, MV-algebras and orthomodular lattices. We will see that this task is far from straightforward. Indeed, Basic algebras can not be represented as semirings since they do not satisfy both distributivity laws, but right-distributivity only; in addition, multiplication need not to be associative in general. 

These observations seem to suggest that a substantial weakening of the concept of semiring would be required to embrace such algebras. An appropriate generalization can be found in \cite{ChajdaL15, ChajdaL} where H. L{\"{a}}nger and one of the present authors discuss the concept of \emph{near semiring}. Taking up ideas from \cite{DiNola13} and \cite{DiNola05}, in order to provide a semiring-like representation of basic algebras, we specialize the concept of near semiring and introduce the notion of \L ukasiewicz near semiring and orthomodular near semiring. \\

\noindent The paper is structured as follows: in section\ \ref{sec: 2} we introduce the notions of near semiring, near semiring with involution and \L ukasiewicz near semiring and discuss some basic properties of these three classes. In section\ \ref{sec: 3} after a concise presentation of basic algebras, we prove that they can be represented by \L ukasiewicz near semirings. In section\ \ref{sec: 4} we discuss several universal algebraic properties of \L ukasiewicz near semirings: congruence regularity, congruence permutability and congruence distributivity. In section\ \ref{sec: 5} we introduce the concept of orthomodular near semiring, and we show that orthomodular lattices can be represented by of these algebraic structures.
Finally, in section\ \ref{sec: 6}, we claim that the variety of involutive integral near semirings is a Church variety \cite{Sal}. This yields an explicit description of central elements and, consequently, a series of direct decomposition theorems.

\section{Near semirings}\label{sec: 2}

\begin{definition}\label{def: near semiring}
A near semiring is an algebra $ \mathbf{R}=\langle R, +, \cdot , 0, 1\rangle $ of type $ \langle 2, 2, 0, 0\rangle $ such that
\begin{itemize}
\item[(i)] $ \langle R, +, 0\rangle  $ is a commutative monoid;
\item[(ii)] $ \langle R, \cdot, 1 \rangle $ is a groupoid satisfying $ x\cdot 1 = x = 1\cdot x $ (a unital groupoid);
\item[(iii)] $ (x+y)\cdot z = (x\cdot z) + (y\cdot z) $;
\item[(iv)] $ x\cdot 0 = 0\cdot x = 0 $.
\end{itemize}
\end{definition}
\noindent
We will refer to the operations $ + $ and $ \cdot $ as sum and multiplication, respectively, and we call the identity in (iii) \emph{right distributivity}. Near semirings generalize semirings to a non-associative and weakly-distributive context. Indeed, a semiring is a near semiring such that $ \langle R, \cdot, 1\rangle $ is a monoid (i.e. $ \cdot $ is also associative) that satisfies \emph{left distributivity}: $ x\cdot (y+z) = (x\cdot y) + (x\cdot z) $, for all $ x, y, z\in R $. Throughout the paper, a near semiring $ \mathbf{R} $ is called \emph{associative} if it satisfies $(x\cdot y)\cdot z=x\cdot( y\cdot z)$, \textit{commutative} if it satisfies $ x\cdot y = y\cdot x $; \textit{idempotent} if it satisfies $ x + x = x $ and \textit{integral} if $ x+ 1 = 1 $ holds.
\begin{remark}\label{rem: order x+y =y}
\emph{Let $ \mathbf{R} $ be an idempotent near semiring. Then $ \langle R, + \rangle $ is a semilattice. In particular, $ \langle R, +\rangle $ can be considered as a join-semilattice, where the induced order is defined as $ x\leq y $ iff $ x+y=y $ and the constant $ 0 $ is the least element. Moreover, whenever $ \mathbf{R} $ is \textit{integral}, the constant 1 is the greatest element with respect to the induced order $ \leq $. }
\end{remark}

\begin{remark}\label{rem: bisemilattices as near semirings}
\emph{Let $ \mathbf{R} $ be an idempotent commutative semiring, whose multiplication is also idempotent ($ x\cdot x = x $). Then, $ \langle R,\cdot \rangle $ is also a semilattice, in particular a meet-semilattice. Notice that, in general, $ \langle R, +, \cdot\rangle $ need not be a lattice. Indeed, the absorption laws may fail.\footnote{Let us remark that, since near semirings satisfy right distributivity only, we may have different forms of absorption.} 
 Moreover, the order induced by the multiplication, $ x\leqslant y $ iff $ x\cdot y = x $, may differ from $ \leq $, as Example \ref{ex: near semiring failing absorption} shows.}
\end{remark}

\begin{example}\label{ex: near semiring failing absorption}
\emph{Let $ \mathbf{R} $ be the near semiring whose universe is $ R=\{0,1,a\} $  and whose sum and multiplication are defined in the following tables:} \\

\vspace{10pt}
\begin{center}
\begin{minipage}{5cm}\begin{tabular}{r|ccc}

$+$ & $0$ & $ a $ & $1$ \\
\hline
$0$ & $0$ & $a$ & $1$ \\
$a$ & $a$ & $a$ & $a$ \\
$1$ & $1$ & $a$ & $1$ \\

\end{tabular}\end{minipage}
\begin{minipage}{5cm}\begin{tabular}{r|ccc}

$ \cdot $ & $0$ & $ $a$ $ & $1$ \\
\hline
$0$ & $0$ & $0$ & $0$ \\
$a$ & $0$ & $a$ & $a$ \\
$1$ & $0$ & $a$ & $1$ \\
\end{tabular}
\end{minipage}
\end{center}
\vspace*{15pt}
\noindent
\emph{It is not difficult to verify that $ \mathbf{R} $ is both additively and multiplicatively idempotent, commutative and associative; therefore, $ \langle R,\cdot \rangle $ is a meet-semilattice. However, $ \mathbf{R} $ is not integral, since $ a + 1=a\neq 1$; and the absorption laws do not hold: $ 1 \cdot (a + 1) = 1\cdot a = a\neq 1 $. It can be seen in Figure \ref{fig:1} that the orders induced by $+$ and $\cdot$ are different.} 
\vspace{15pt}

\begin{figure}[h]
\begin{tikzpicture}\label{fig:1}

\draw (0,0) -- (0,1.5); 
\draw (0,1.5) -- (0,3); 
\draw (0,0) node {$\bullet$};
\draw (-0.5, 0) node {0};
\draw (0,1.5) node {$ \bullet$};
\draw (-0.5, 1.5) node {1};
\draw (0.6, 1.1) node {$\leq$};
\draw (0,3) node {$\bullet$};
\draw (-0.5, 3) node {a}; 

\draw (4,0) -- (4,1.5); 
\draw (4,1.5) -- (4,3); 
\draw (4,0) node {$\bullet$};
\draw (4.5,0) node {0};
\draw (4,1.5) node {$ \bullet$};
\draw (4.5,1.5) node {a};
\draw (3.4, 1.1) node {$\leqslant$};
\draw (4,3) node {$\bullet$};
\draw (4.5, 3) node {1};

\end{tikzpicture}
\caption{The Hasse diagrams of the partial orders induced by sum, $ \leq $ (left hand side), and multiplication, $ \leqslant $ (right hand side).}
\end{figure}
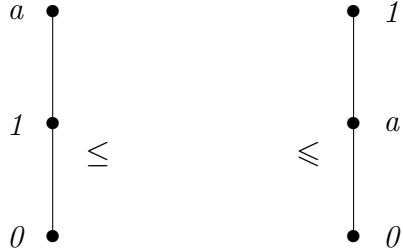
\end{example}

The following lemma states that in any near semiring, multiplication is monotone on the right hand side, due to right distributivity.
\begin{lemma}\label{lemma: monotonicity of multiplication}
Let $ \mathbf{R} $ be a near semiring. Then, $ x\leq y $ implies $ x\cdot z\leq y\cdot z $. 
\begin{proof}
Suppose $ x\leq y $, i.e. $ x+ y = y $. Therefore $ y\cdot z = (x + y)\cdot z = (x\cdot z) + (y\cdot z) $, which implies that $ x\cdot z\leq y\cdot z $. 
\end{proof}
\end{lemma}

Since near semirings are, in general, not distributive (left distributivity does not hold), multiplication is not monotone in the left component. However, we there will be cases in which distributivity holds for some special elements.

\begin{definition}\label{def: near semi with involution}
Let $ \langle R,+,\cdot, 0,1\rangle $ be an \emph{idempotent} near semiring, with $ \leq $ the induced order. A map $ \alpha: R\rightarrow R $ is called an \emph{involution} on R if it satisfies the following conditions for each $ x,y \in R $:
\begin{itemize}
\item[(a)] $ \alpha(\alpha (x)) = x $; 
\item[(b)] if $ x\leq y $ then $ \alpha(y)\leq\alpha (x) $.
\end{itemize}
The algebra $ \mathbf{R} =\langle R, +,\cdot , 0,1,\alpha\rangle $ will be called an \emph{involutive near semiring}.
\end{definition}
\noindent
Sometimes, if no confusion is possible, we will write $\alpha\alpha x$ in place of $\alpha(\alpha (x))$.
Some basic arithmetical properties of involutive near semirings are presented in the following lemma.
\begin{lemma}\label{lemma: aritmetica}
Let $ \mathbf{R} $ be an involutive near semiring . Then 

\begin{itemize}
\item[(i)] $ \alpha (x + y) + \alpha (x) = \alpha (x) $.
\item[(ii)] $ \mathbf{R} $ is \emph{integral} if and only if $ \alpha(0) = 1 $ (and consequently $ \alpha(1) = 0 $).
\end{itemize}
\begin{proof}
(i) Since + is idempotent, we have that $ \alpha(x) = \alpha (x) + \alpha (x) $. Moreover, since $x + y=(x+x)+y=x+(x+y)$, $ x\leq x + y $ and $ \alpha $ is an involution on $ R $, $ \alpha(x + y)\leq\alpha (x) $. Therefore $ \alpha (x + y) +\alpha (x) \leq \alpha(x) + \alpha (x) = \alpha (x) $. The converse $ \alpha (x) \leq \alpha (x) + \alpha (x + y) $ holds because $ \langle R, +\rangle $ is a join-semilattice as we noticed in Remark \ref{rem: order x+y =y}. \\
(ii) Suppose that $ \mathbf{R} $ is integral, i.e. $ x\leq 1 $ for each $ x\in R $, then $ \alpha (1)\leq \alpha (x) $. Since $ \alpha $ is an involution we have that $ \alpha(1) + x = x $ for each $ x\in R $, which means that $ \alpha (1) $ is a neutral element with respect to the sum and since $ \langle R, +, 0\rangle $ is a (commutative) monoid, the neutral is unique\footnote{This is in fact true for semigroups.}, thus $ \alpha (1) = 0 $ and $ \alpha (0) = 1 $. For the converse, suppose $ \alpha (0) = 1 $. Then, by (i), we have that $ \alpha (x + y) \leq \alpha (x) $, which, for $ x = 0 $, implies $ \alpha (y)\leq \alpha(0) = 1 $, which means that for each $ x\in R $ we have $ x\leq 1 $, i.e. $ \mathbf{R} $ is integral.
\end{proof}
\end{lemma}

\begin{remark}\label{rem: 1 top non accade sempre}
\emph{Notice that, in general, $ x\leq\alpha (0) $. However, it is not true that in any involutive near semiring  $\alpha(0)=1$ (see Example \ref{ex: not-int}). In fact, there are cases of non integral involutive near semiring}
\end{remark}

\begin{example}\label{ex: not-int}\text{}\\
\begin{center}
 \begin{tabular}{r|rrr}
$'$ & $0$ & $1$ & $2$\\
\hline
   & $2$ & $1$ & $0$
\end{tabular} \hspace{.5cm}
\begin{tabular}{r|rrr}
$+$ & $0$ & $1$ & $2$\\
\hline
    $0$ & $0$ & $1$ & $2$ \\
    $1$ & $1$ & $1$ & $2$ \\
    $2$ & $2$ & $2$ & $2$
\end{tabular} \hspace{.5cm}
\begin{tabular}{r|rrr}
$\cdot$ & $0$ & $1$ & $2$\\
\hline
    $0$ & $0$ & $0$ & $0$ \\
    $1$ & $0$ & $1$ & $2$ \\
    $2$ & $0$ & $2$ & $2$
\end{tabular} 
\end{center}
\end{example}

%

\begin{theorem}\label{th: dual near semiring}
Let $ \mathbf{R} $ be an involutive near semiring and define two new operations as $ x +_{\alpha} y = \alpha (\alpha (x) + \alpha (y)) $ and $ x\cdot_{\alpha} y= \alpha (\alpha(x)\cdot\alpha(y)) $. Then:
\begin{itemize}
\item[(a)] $ x+ y= \alpha(\alpha(x)+_{\alpha} \alpha(y)) $, $ x\cdot y= \alpha(\alpha(x)\cdot_{\alpha} \alpha(y)) $;
\item[(b)] $\mathbf{R}_{\alpha}= \langle R, +_{\alpha}, \cdot_{\alpha}, \alpha, \alpha(0), \alpha(1)\rangle $ is an involutive near semiring.
\end{itemize}
\begin{proof}
(a) By definition of $ +_{\alpha} $ we have that $\alpha(\alpha(x)+_{\alpha} \alpha(y))=\alpha\alpha(\alpha\alpha(x)+ \alpha\alpha(y))=x+y$. The proof runs analogously for $ \cdot_{\alpha} $. \\
(b) We start by showing that $ \langle R, +_{\alpha}, \alpha (0)\rangle $ is a commutative monoid. Commutativity of $ +_{\alpha} $ trivially follows by definition. Furthermore:
\begin{align*}
& (x+_{\alpha} y) +_{\alpha} z = \alpha(\alpha(x +_{\alpha} y) + \alpha(z)) && {{(\text{Def.}}\; +_{\alpha})}
\\ & = \alpha(\alpha\alpha(\alpha(x) + \alpha(y)) + \alpha(z)) && ({\text {Def.} }\; +_{\alpha})
\\ & = \alpha ((\alpha (x) + \alpha (y)) + \alpha(z)) && (\text{Inv.}) 
\\ & = \alpha (\alpha(x) + (\alpha(y) + \alpha(z))) && (\text{Ass.}\; +)
\\ & = \alpha(\alpha (x) + \alpha(y +_{\alpha} z)) && (\text{Def.}\; +_{\alpha}) 
\\ & = x +_{\alpha} (y +_{\alpha} z) && (\text{Def} \; +_{\alpha}),
\end{align*}  
proving associativity of $ +_{\alpha} $. Finally, 
\begin{align*}
& x +_{\alpha} \alpha (0) = \alpha(\alpha(x) + \alpha(\alpha(0))) && (\text{Def.} \; +_{\alpha})
\\ & = \alpha (\alpha(x) + 0) && (\text{Inv.}) 
\\ & = \alpha(\alpha(x)) = x && (\text{Monoid})
\end{align*} 
The fact that $ \alpha(0) $ is also a left neutral follows from commutativity. The proof of the fact that $ \langle R, \cdot_{\alpha},\alpha(1)\rangle $ is a groupoid with $ \alpha(1) $ as neutral element is analogous.  \\
Furthermore, $ x\cdot_{\alpha} \alpha(0) = \alpha (\alpha (x) \cdot \alpha(\alpha (0)) = \alpha (\alpha (x) \cdot 0) = \alpha (0) $ and similarly to show that $ \alpha(0)\cdot_{\alpha} x = \alpha(0) $. \\
It only remains to show that right distributivity holds. 
\begin{align*}
& (x +_{\alpha} y)\cdot_{\alpha} z = \alpha (\alpha (x+_{\alpha} y)\cdot \alpha(z)) && (\text{Def.} \; \cdot_{\alpha}) 
\\ & = \alpha(\alpha\alpha(\alpha(x)+\alpha(y))\cdot \alpha (z)) && (\text{Def.} \; +_{\alpha})
\\ & = \alpha ( (\alpha (x) + \alpha (y)) \cdot \alpha (z)) && (\text{Inv.}) 
\\ & = \alpha ((\alpha (x) \cdot \alpha (z)) + (\alpha (y)\cdot\alpha(z))) && (\text{Distr.})
\\ & = (x\cdot_{\alpha} z) +_{\alpha} (y \cdot_{\alpha} z) && (\text{Def.}) 
\end{align*}
Therefore $ \langle R, +_{\alpha}, \cdot_{\alpha}, \alpha(0), \alpha(1)\rangle $ is a near semiring.
\end{proof}
\end{theorem}

In general, for a given near semiring $ \mathbf{R} $, we will refer to $ \mathbf{R}_{\alpha} $ as the \emph{dual} near semiring. 

As we mentioned in the introduction, semiring-like structures are relevant to the theory of prominent algebraic structures from many-valued logics. In particular, it was shown by Belluce, Di Nola, Ferraioli \cite{DiNola13} and Gerla \cite{Gerla03} that MV-algebras can be represented as semirings. Following the same idea, we aim at representing some algebraic structures deriving from quantum logics as near semirings. For this reason we introduce the notion of \L ukasiewicz near semiring.
\begin{definition}\label{def: Luka near semiring}
Let $ \mathbf{R} $ be an involutive near semiring. $ \mathbf{R} $ is called a \emph{\L ukasiewicz near semiring} if it satisfies the following additional identity
\begin{itemize}
\item[(\L)] $ \alpha (x\cdot \alpha (y))\cdot \alpha (y) = \alpha (y\cdot \alpha(x))\cdot \alpha (x) $.
\end{itemize}
A \emph{semiring} satisfying \emph{(\L)} will be called a
 \emph{\L ukasiewicz semiring}.
\end{definition}
\noindent
Identity (\L), in Definition \ref{def: Luka near semiring}, clearly reflects \L ukasiewicz identity in the standard axiomatization of MV-algebras. As already noticed, the constant $ 1 $ need not necessarily be the top element with respect to the order $ \leq $ in general. However, as the next lemma shows, this is always the case for \L ukasiewicz near semiring. This fact will be frequently used throughout this paper.

\begin{lemma}\label{lemma: Lukas near semi sono integrali}
Let $ \mathbf{R} $ be \L ukasiewicz near semiring. Then
\begin{itemize}
\item[(a)] $ x\cdot \alpha(x) = \alpha (x)\cdot x = 0 $;
\item[(b)] $ \mathbf{R} $ is integral;
\item[(c)] $ x\cdot \alpha (x + y) = 0 $;
\item[(d)] $ (x+y)\cdot \alpha (x) =  y\cdot \alpha(x) $;
\item[(e)] $ x+y=\alpha (\alpha(x\cdot \alpha(y))\cdot \alpha(y)) $. 
\end{itemize}
\begin{proof}
(a) Let us observe that, upon setting $x = 0 $ and $ y=1$ in (\L), we get $ \alpha (0)\cdot \alpha (1) = \alpha(0\cdot \alpha(1))\cdot \alpha(1) = \alpha (1\cdot \alpha (0))\cdot \alpha (0) = \alpha (\alpha (0))\cdot \alpha (0) = 0 $.  Since $ 0 $ is the unit with respect to the sum, we have that $ 0 +\alpha(x) = \alpha (x) $, i.e. $ 0\leq \alpha (x) $. Therefore $ x\leq \alpha (0) $ and then $ x + \alpha (0) = \alpha (0) $. Using these two facts, we obtain
\begin{align*}
& 0 = (x + \alpha (0))\cdot \alpha(1) 
\\ & = (x\cdot \alpha (1)) + (\alpha(0)\cdot \alpha (1)) && (\text{Distr.})
\\ & = (x\cdot \alpha (1)) + 0 
\\ & = x\cdot \alpha (1)
\end{align*}
We finally get 
\begin{align*}
& x \cdot \alpha (x) = \alpha (\alpha (x))\cdot \alpha (x) 
\\ & = \alpha(1\cdot \alpha (x))\cdot \alpha (x) 
\\ & = \alpha (x\cdot \alpha (1))\cdot \alpha (1) && (\text{\L})
\\ & = \alpha (0)\cdot \alpha (1) 
\\ & = 0.
\end{align*}
\noindent
 Since $\alpha$ is an involution, it follows that also $ \alpha(x)\cdot x = 0 $. \\
(b) $ \alpha(1) = 1\cdot \alpha (1) = 0 $, by (a). Then by Lemma \ref{lemma: aritmetica} (ii) we have that $ \mathbf{R} $ is integral. \\
(c) Since $ x\leq x+ y $. Then, by Lemma \ref{lemma: monotonicity of multiplication}, $ x\cdot \alpha (x +y)\leq (x+y)\cdot \alpha (x + y) = 0 $. \\
(d) It is enough using right distributivity and (a), indeed $ (x+y)\cdot \alpha (x) = (x\cdot \alpha(x)) + (y\cdot\alpha (x)) = 0 + (y\cdot\alpha (x)) =  y\cdot \alpha(x) $. \\
(e) 
\begin{align*}
& \alpha (x\cdot\alpha (y))\cdot\alpha(y) = \alpha (y\cdot \alpha (x))\cdot\alpha (x) && (\text{\L})
\\ & = \alpha ((x+y)\cdot \alpha(x))\cdot\alpha (x) && (\text{Item (d)}) 
\\ & = \alpha(x \cdot \alpha (x+y))\cdot \alpha(x + y) && (\text{\L})
\\ & = \alpha (0) \cdot \alpha (x+y) && (\text{Item (c)}) 
\\ & = 1\cdot \alpha(x+y) && (\text{Item (b)}) 
\\ & = \alpha(x+y).
\end{align*}
Therefore $ x + y= \alpha (\alpha(x\cdot \alpha(y))\cdot \alpha(y)) $
\end{proof} 
\end{lemma}
\noindent
%
%
The next lemma shows that, in the specific case of \L ukasiewicz near semirings, the order induced by the sum is equivalently expressed by multiplication.
\begin{lemma}\label{lem: x leq y iff x * y' = 0}
Let $ \mathbf{R} $ be a \L ukasiewicz near semiring. Then $ x\leq y $ if and only if $ x\cdot\alpha(y) = 0 $. 
\begin{proof}
Let $ a\leq b $, for some $a, b\in R $. Then $ a + b = b $ and, by Lemma \ref{lemma: Lukas near semi sono integrali}-(c), we get that $ 0 = a\cdot \alpha (a + b) = a\cdot \alpha (b) $. \\
Conversely, suppose that $ a\cdot \alpha(b) = 0 $ for some $ a,b\in R $. 
\begin{align*}
& a + b = \alpha (\alpha(a\cdot\alpha(b))\cdot \alpha (b))) && (\text{Lemma \ref{lemma: Lukas near semi sono integrali}(e)})
\\ & = \alpha (\alpha (0)\cdot \alpha (b)) && (\text{Assumption})
\\ & = b. 
\end{align*}
Therefore $ a\leq b $. 
\end{proof}
\end{lemma}
\noindent
Theorem \ref{th: dai near semirings ai semirings} clarifies the connection between \L ukasiewicz near semirings and \L ukasiewicz semirings.

\begin{theorem}\label{th: dai near semirings ai semirings}
Let $\mathbf{ R} $ be a \L ukasiewicz near semiring whose multiplication is associative. Then multiplication is also commutative, and therefore $ \mathbf{R} $ is a commutative \L ukasiewicz semiring. 
\begin{proof}
Suppose $ \langle R, \cdot\rangle $ is a semigroup. Then 
\begin{align*}
& \alpha(x\cdot y)\cdot (y\cdot x) = (\alpha (x\cdot y)\cdot y)\cdot x && (\text{Assumption})
\\ &= (\alpha (\alpha(y)\cdot \alpha(x))\cdot \alpha (x)) \cdot x && (\text{\L}) 
\\ & = (  \alpha (\alpha(y)\cdot \alpha(x))\cdot (\alpha (x) \cdot x) && (\text{Assumption}) 
\\ & =  (  \alpha (\alpha(y)\cdot \alpha(x))\cdot 0 && (\text{Lemma \ref{lemma: Lukas near semi sono integrali}}) 
\\ 0.
\end{align*}
Therefore $ \alpha(x\cdot y)\cdot (y\cdot x) = 0 $. Analogously, $ \alpha(y\cdot x)\cdot (x\cdot y) = 0 $. 
Applying Lemma \ref{lem: x leq y iff x * y' = 0} to both the equations, we obtain that $ \alpha (x\cdot y)\leq\alpha (y\cdot x) $ and $ \alpha (y\cdot x)\leq\alpha(x\cdot y) $. Therefore $ \alpha(x\cdot y) = \alpha (y\cdot x) $, i.e. $ x\cdot y = y\cdot x $. 
Therefore, multiplication commutes. Hence, to prove that $ \mathbf{R} $ is a \L ukasiewicz semiring, it suffices to observe that left distributivity follows straight away from right distributivity.
\end{proof}
\end{theorem}
\noindent
As an immediate consequences of the previous result we obtain that
\begin{corollary}\label{cor: Luka semirings are commutative}
Every \L ukasiewicz semiring is commutative.
\end{corollary} 

\begin{corollary}\label{cor: semirings 2}
A \L ukasiewicz near semiring is a \L ukasiewicz semiring if and only if multiplication is associative. 
\end{corollary}
\noindent
By Lemma \ref{lemma: Lukas near semi sono integrali}, any \L ukasiewicz near semiring is integral. Therefore, it makes sense to introduce the notion of {\it interval} on a \L ukasiewicz near semiring $ \mathbf{R} $: $ [a,1]=\{ x\in R | a\leq x \}$. The next result shows that any such interval can be equipped with an antitone involution.

\begin{theorem}\label{th: antitone involution on intervals for Lukas. semirings}
Let ${\mb{R}} $ be a \L ukasiewicz near semiring, $ \leq $ the induced order, and $ a\in R $. The map $ h_{a}: [a,1]\rightarrow [a,1] $, defined by $ x\mapsto x^{a}= \alpha( x\cdot\alpha(a)) $ is an antitone involution on the interval $ [a,1] $. 
\begin{proof}
We first show that $ h_{a} $ is well defined. Indeed, since $ \mathbf{R} $ is integral (Lemma \ref{lemma: Lukas near semi sono integrali}) we have that $ x\leq 1 $, thus $ x\cdot\alpha (a)\leq 1\cdot\alpha (a) $ by monotonicity, then $ a= \alpha(\alpha (a)) =\alpha (1\cdot \alpha (a))\leq \alpha (x\cdot\alpha (a))= x^{a} $, i.e. $ x^{a}\in [a,1] $. 
Moreover, $ h_{a} $ is antitone. Suppose $ x,y\in [a,1] $ with $ x\leq y $. Since multiplication is monotone (Lemma \ref{lemma: monotonicity of multiplication}) we get that $ x\cdot \alpha (a) \leq y\cdot\alpha (a) $. Therefore $ y^{a}= \alpha (y\cdot\alpha (a)) \leq \alpha (x\cdot\alpha (a)) = x^{a} $, i.e. $ h_{a} $ is antitone. \\
Since for any $ x\in [a,1] $, $ a\leq x $, i.e. $ a + x = x $, then by Lemma \ref{lemma: Lukas near semi sono integrali}-(c)
\begin{align*}
& a\cdot \alpha (x) = a\cdot \alpha (a+x) = 0 && (*)
\end{align*} 
From this fact we obtain that:
\begin{align*}
& x^{aa}= \alpha (x^{a}\cdot \alpha (a)) = \alpha (\alpha (x\cdot \alpha(a))\cdot\alpha (a)) && (\text{Definition}) 
\\ & = \alpha(\alpha (a\cdot\alpha (x))\cdot\alpha (x)) && (\text{\L}) 
\\ & = \alpha (\alpha(0)\cdot\alpha (x) )&& (*) 
\\ & = \alpha (1\cdot \alpha (x)) && (\text{Integrality}) 
\\ & = \alpha(\alpha (x)) = x.
\end{align*}
This shows that $ h_{a} $ is an antitone involution on the interval $ [a,1] $.
\end{proof}
\end{theorem}
\noindent
Let us observe that, in \cite{Chajda15, Chajda09interval}, the involution constructed in Theorem \ref{th: antitone involution on intervals for Lukas. semirings} is termed \emph{sectional involution}.

\section{Basic algebras as near semirings}\label{sec: 3}
Basic algebras, introduced in the last decade by Hala{\v{s}}, K{\"u}hr, and one of the authors of this paper, as a common generalization of both MV-algebras and orthomodular lattices. They can be regarded as a non-associative and non-commutative generalization of MV algebras. These algebras are in bijective correspondence with bounded lattices having an antitone involution on every principal filter (\emph{sectional antitone involutions}). An introductory as well as comprehensive survey on basic algebras can be found in \cite{Chajda15}.

In this section we discuss the relations between \L ukasiewicz near semirings and basic algebras. Let us recall that a basic algebra is an algebra $ \mathbf{A}=\langle A, \oplus, ', 0\rangle $ satisfying the following identities: 
\begin{itemize}
\item[(BA1)] $ x \oplus 0 = x $;
\item[(BA2)] $ x'' = x $;
\item[(BA3)] $ (x'\oplus y)' \oplus y = (y \oplus x')' \oplus x $;
\item[(BA4)] $ (((x\oplus y)' \oplus y)' \oplus z) ' \oplus (x\oplus z) = 1 $.
\end{itemize}
where $ 0' = 1 $ and (BA3) is the \textit{\L ukasiewicz identity}. It is not difficult to show that every MV algebra is a basic algebra. More precisely the class of MV algebras is a subvariety of the variety of basic algebras and it is axiomatized by the identity expressing associativity of $ \oplus $,\footnote{Indeed, it is shown in \cite{Chajda15} that if $ \oplus $ is associative then it is also commutative.} (for details see \cite{Chajda15}). Every basic algebra is in fact a bounded lattice, where the lattice order is defined by $ x\leq y $ iff $ x'\oplus y = 1 $, the join operation is defined by $ x\vee y =  (x' \oplus y)'\oplus y $, while the meet is defined \'a la de Morgan:  $ x\wedge y= (x' \vee y')' $. It can be verified that $ 0 $ and $ 1 $ are the bottom and the top elements, respectively, of the lattice. 

Conversely, let us remark that, in any bounded lattice with sectional antitone involutions $ \langle L, \vee , \wedge , (^{a})_{a\in L} ,  0 ,1 \rangle $ (the interested reader may consult \cite{Chajda07}, \cite{Chajda15} ), it is possible to define the operations
\begin{equation}\label{eq: operazioni di basic algebras da un bounded lattice}
  x' = x^{0}, \;\;\;  x\oplus y:= (x^{0}\vee y)^{y}, 
\end{equation}  
such that $ \langle L, \oplus, ', 0,1\rangle $ is a basic algebra. It can be proven that this correspondence is one to one.
We will use this fact to establish a correspondence between \L ukasiewicz near semirings and basic algebras. In particular, we will discuss how basic algebras can be represented in terms of near semirings.
\begin{theorem}\label{th: basic algebras dai Luk near semirings}
If $ \mathbf{R} $ is a \L ukasiewicz near semiring, then the structure $ \mathbf{B(R)}= \langle R, \oplus, \alpha , 0\rangle $, where $ x\oplus y$ is defined by $\alpha((\alpha (x) + y)\cdot \alpha (y))$ is a basic algebra. 
\begin{proof}
The reduct $ \langle R, +, 1\rangle $ is a (join) semilattice whose top element is $1$ (Remark \ref{rem: order x+y =y} and Lemma \ref{lemma: Lukas near semi sono integrali}-(b)). From Theorem \ref{th: dual near semiring}, we have that $ \langle R, +_{\alpha}\rangle $ is the dual meet semilattice. Therefore $ \langle R, +, +_{\alpha}, 0,1\rangle $ is a bounded lattice. Furthermore, by Theorem \ref{th: antitone involution on intervals for Lukas. semirings} the map $ x\mapsto x^{a}= \alpha (x\cdot\alpha (a)) $ is an antitone involution on the interval $ [a,1] $ for all $ a\in R $. So $ \langle R, +, +_{\alpha}, (^{a})_{a\in R} , 0,1\rangle $ is a bounded involution lattice with sectional antitone involutions. And therefore it can be made into a basic algebra upon setting the operations as in equations \eqref{eq: operazioni di basic algebras da un bounded lattice}. It follows that $ x^{0}= x' = \alpha(x) $ and $ x\oplus y = \alpha ((\alpha(x) + y)\cdot \alpha (y)) $.
\end{proof}
\end{theorem}
\noindent
The next result shows a converse of the previous theorem: any basic algebra induces a \L ukasiewicz near semiring.

\begin{theorem}\label{th: da basic algebras a Lukas near semiring}
If $ \mathbf{B}= \langle B,\oplus, ', 0\rangle $ is a basic algebra, then the structure $ \mathbf{R}(\mathbf{B})=\langle B, + , \cdot , \alpha , 0, 1\rangle $, where $x+y$, $ x\cdot y$ and $\alpha (x)$ are defined by $(x' \oplus y)'\oplus y$, $(x'\oplus y')'$, $x'$, and $1=0'$, respectively, is a \L ukasiewicz near semiring.
\begin{proof}
As we mentioned, in any basic algebras $ (x' \oplus y)'\oplus y $ defines a join-semilattice, whose least and greatest elements are, respectively, 0 and 1. This assures that $ \langle R, + , 0\rangle $ is a commutative monoid. Furthermore, it was shown in \cite{Chajda15} that $ (x + y)\cdot z =  (x\cdot z) + (y\cdot z) $. Let us now prove that $1$ is a neutral element for the multiplication. We note that 
\begin{align*}
& x\cdot 1 = (x' \oplus 1')' && (\text{Def.} \; \cdot )
\\ & = (x'\oplus 0)' && (\text{Int.})
\\ & = x'' && (\text{BA1}) 
\\ & = x && (\text{BA2})
\end{align*}
Upon observing that, in basic algebras, $ x\oplus 0 = x $ (BA1) implies $ 0\oplus x = x $, then one analogously proves that $ 1\cdot x = x $. To prove that 0 is an annhilator of multiplication, we show that
\begin{align*}
& x\cdot 0 = (x' \oplus 0')' && (\text{Def.} \; \cdot ) 
\\ & =(x'\oplus 1) ' && (\text{Int.}) 
\\ & =1' && (\text{BA}) 
\\ & = 0.
\end{align*}
The proof that $ 0\cdot x = 0 $ is analogous. Therefore $ \mathbf{R}(\mathbf{B}) $ is a near semiring. Since $ \alpha (x) = x' $, it is clear that it is also an antitone involution. We are left with checking that $ \mathbf{R}(\mathbf{B}) $ satisfies the conditions of a \L ukasiewicz near semiring, Definition \ref{def: Luka near semiring}. \\
As regards condition (\L),
\begin{align*}
& \alpha(x\cdot\alpha (y))\cdot\alpha(y) = ((x\oplus y)'\oplus y)' && (\text{Def.})
\\ & = ((y\oplus x)'\oplus x)' && (\text{BA3}) 
\\ & =  \alpha(y\cdot\alpha (x))\cdot\alpha(x) && 
\end{align*}  
This concludes the proof that $ \mathbf{R}(\mathbf{B}) $ is a \L ukasiewicz near semiring. 
\end{proof}
\end{theorem} 
\noindent
\noindent The results above state a correspondence between near \L ukasiewicz semirings and basic algebras.
In order to analyze these maps properly, we will refer to the variety of basic algebras and of \L ukasiewicz near semiring as $ \mathcal{B} $ and $ \mathcal{R} $, respectively. In Theorem \ref{th: da basic algebras a Lukas near semiring} we considered a map $ f:\mathcal{B}\rightarrow\mathcal{R} $ associating to each basic algebra a \L ukasiewicz near semiring $ \mathbf{R}(\mathbf{B}) $. On the other hand, in Theorem \ref{th: basic algebras dai Luk near semirings}, we applied a map $ g:\mathcal{R}\rightarrow\mathcal{B} $, associating to any \L ukasiewicz near semiring $ \mathbf{R} $ a basic algebra $ \mathbf{B}(\mathbf{R}) $.

The next theorem shows that $ \mathbf{B}(\mathbf{R}(\mathbf{B})) $ actually coincides with $ \mathbf{B} $ and, viceversa, $ \mathbf{R} $ coincides with  $ \mathbf{R}(\mathbf{B}(\mathbf{R})) $. 
\begin{theorem}\label{th: f e g mutually inverse}
The maps $ f $ and $ g $ are mutally inverse.
\begin{proof}
We start by checking that $ \mathbf{B}(\mathbf{R}(\mathbf{B}))=\mathbf{B} $. 
We first note that $ f(x')=\alpha (x) $ and $ g(\alpha(x)) = x' $. Therefore $ f(g(\alpha(x))) = \alpha (x) $ and $ g(f(x')) = x' $. We have to prove that $ x\oplus y = x\;\widehat{\oplus}\; y $, where by $ x\;\widehat{\oplus}\; y $ we indicate the sum in $ \mathbf{B}(\mathbf{R}(\mathbf{B})) $.  We use the fact that $ \mathbf{R}(\mathbf{B}) $ is a \L ukasiewicz near semiring (Theorem \ref{th: da basic algebras a Lukas near semiring}), whose sum and multiplication are indicated by $ \widehat{+} $ and $ \widehat{\cdot} $, respectively. 
\begin{align*}
& x\;\widehat{\oplus}\; y = \alpha((\alpha(x)\; \widehat{+}\; y)\; \widehat{\cdot}\;\alpha(y)) && (\text{Def.})
\\ & = \alpha((\alpha(x)\;\widehat{\cdot}\;\alpha(y))\;\widehat{+}\;(y\;\widehat{\cdot}\;\alpha(y))) && (\text{Distr.})
\\ &  = \alpha((\alpha(x)\;\widehat{\cdot}\;\alpha(y))\;\widehat{+}\; 0) && (\text{Lemma \ref{lemma: Lukas near semi sono integrali}}) 
\\ & =\alpha((\alpha(x)\;\widehat{\cdot}\;\alpha(y)) = x\oplus y.
\end{align*}
This is enough to have that $ \mathbf{B}(\mathbf{R}(\mathbf{B}))=\mathbf{B} $. To see that $ \mathbf{R}(\mathbf{B}(\mathbf{R})) =\mathbf{R} $ we need to check that $ x\;\widehat{+}\; y = x + y $ and $ x\;\widehat{\cdot}\; y = x\cdot y $. We begin with the latter equality: $ x\;\widehat{\cdot}\; y = (x'\;\widehat{\oplus}\;y')' = x\cdot y $. Concerning the former, we have that $ x\;\widehat{+}\; y =(x'\;\widehat{\oplus}\;y)'\;\widehat{\oplus}\; y = \alpha((\alpha(x\cdot\alpha(y)))\cdot\alpha(y)) = x + y $ by Lemma \ref{lemma: Lukas near semi sono integrali}-(e). 
\end{proof}
\end{theorem}

As a corollary of the representation of basic algebras as \L ukasiewicz near semirings, we get the one-to-one correspondence between MV-algebras and the variety of commutative \L ukasiewicz near semiring.
The following results readily follow from Theorems \ref{th: basic algebras dai Luk near semirings} and \ref{th: da basic algebras a Lukas near semiring} and the fact that a basic algebra is an MV-algebra if and only if $ \oplus $ is associative.
\begin{corollary}\label{cor: MV and near Lukas semirings}
Let $ \mathbf{M}=\langle M, \oplus , ', 0\rangle $ an MV-algebra. Then the structure $ \mathbf{R}(\mathbf{M})=\langle B, + , \cdot , \alpha , 0, 1\rangle $, where $x+y$, $ x\cdot y$ and $\alpha (x)$ are defined by $(x' \oplus y)'\oplus y$, $(x'\oplus y')'$, $x'$, and $1=0'$ respectively, is a  \L ukasiewicz semiring.
\end{corollary}
\begin{corollary}\label{cor: MV and near Lukas semirings 2}
Let $ \mathbf{R}=\langle R, +, \cdot , \alpha , 0, 1\rangle $ be a \L ukasiewicz semiring and let $ x\oplus y = \alpha((\alpha(x)+ y)\cdot\alpha(y)) $. Then $ \mathbf{M}(\mathbf{R})=\langle R, \oplus, \alpha , 0\rangle $ is an MV-algebra. 
\end{corollary}
\noindent
Corollaries above are tightly related to a very similar result in \cite{DiNola13}, where it is shown that to every MV-algebra corresponds an MV-semiring: a commutative semiring with involution that satisfies the identity (c) in Lemma \ref{lemma: Lukas near semi sono integrali} (in our terminology) and
\begin{equation}\label{eq: MV-semiring}
x + y = \alpha (\alpha (x)\cdot\alpha (\alpha(x)\cdot y)). 
\end{equation}
\noindent

\section{Congruence Properties of \L ukasiewicz near semirings}\label{sec: 4}

In this section we discuss several congruence properties of \L ukasiewicz near semirings. Recall that an algebra $ \mathbf{A} $ is \emph{congruence regular} if any congruence $ \theta\in Con(\mathbf{A}) $ is determined by any of its cosets; namely if $ \theta, \phi \in Con (\mathbf{A}) $ and $ a\in A $ then 
$$ [a]_{_\theta}=[a]_{_\phi}\;\;\;\text{implies}\;\;\; \theta =\phi $$
A variety $ \mathcal{V} $ is \emph{congruence regular} if every member of $ \mathcal{V} $ is congruence regular. 
 A theorem due to Cs\'ak\'any shows (see \cite{Chajda03book} for details) that a variety $ \mathcal{V} $ is congruence regular if and only if there exists a set of ternary terms $ t_{i}(x,y,z) $ with $ i \geq 1 $ such that  
$$ t_{i}(x,y,z) = z \;\;\text{for any}\; i \;\;\text{if and only if} \;\; x = y $$ 
An algebra $ \mathbf{A} $ is said to be \textit{congruence permutable} if for any two congruences $ \theta, \phi \in Con(\mathbf{A}) $ it holds that $ \theta \circ\phi =\phi\circ\theta $. \\
An algebra $ \mathbf{A} $ is \emph{congruence distributive} if the complete lattice of its congruences is distributive. \\
A variety $ \mathcal{V} $ is congruence permutable (congruence distributive, resp.) if every member of $ \mathcal{V} $ is congruence permutable (congruence distributive, resp.). \\
Finally, an algebra $ \mathbf{A} $ is \emph{arithmetical} if it is both congruence permutable and congruence distributive. A variety $ \mathcal{V} $ is arithmetical if each algebra $ \mathbf{A}\in\mathcal{V} $ is arithmetical. 

It was proven by Mal'cev \cite{Malcev} that congruence permutability is equivalent to the existence of a certain (uniformly defined) term operation. Precisely, a variety $ \mathcal{V} $ is congruence permutable if and only if there exists a ternary term operation $ p(x,y,z) $ such that the identities 
\[
p(x,x,y) = y \;\;\text{and}\;\;\; p(x,y,y) = x
\]
hold in $ \mathcal{V} $. The term $ p $ is usually referred to as a Mal'cev term for $ \mathcal{V} $. \\
Similarly, congruence distributivity is witnessed by the existence of the so-called J\'onsson terms. In particular, a variety $ \mathcal{V} $ is congruence distributive if there exists a ternary term operation $ M(x,y,z) $, for which the identities 
\[
M(x, x, y) = M(x, y, x) = M(y, x, x) = x
\]
hold in $ \mathcal{V} $. $ M $ is usually referred to as a majority term for $ \mathcal{V} $. 

\begin{theorem}\label{th: near Luk semirings are congruence regular}
The variety of \L ukasiewicz near semirings is congruence regular, with witness terms:
\begin{eqnarray*}
t_1(x,y,z) =& ((x \cdot \alpha(y)) + (y\cdot \alpha(x)) + z; \\
t_2(x,y,z) = & \alpha ((x\cdot \alpha (y) ) + (y\cdot \alpha(x))) \cdot z.
\end{eqnarray*}

\begin{proof}
All we need to check is that $ t_1(x,y,z) = t_{2}(x,y,z) = z $ if and only if $ x = y $. Suppose that $ x = y $; then $ t_1 (x,x,z)= ((x \cdot \alpha(x)) + (x\cdot \alpha(x)) + z   =  (0+0)+ z = 0 + z = z $. On the other hand $  t_2(x,x,z) = \alpha ((x\cdot \alpha (x) ) + (x\cdot \alpha(x))) \cdot z  = \alpha(0 + 0) \cdot z = 1\cdot z = z.  $ For the converse, suppose $ t_1(x,y,z) = t_2 (x,y,z) = z $, which, setting $ a = (x \cdot \alpha(y)) + (y\cdot \alpha(x)) $, reads
\begin{align}
&  a + z = z && \label{un}
\\ & \alpha (a) \cdot z = z && \label{du}
\end{align}
Equation \eqref{un} above implies that $ a\leq z $, hence $ \alpha (z)\leq\alpha (a) $. We now claim that $ a = 0 $. Indeed 
\begin{align*}
& a = \alpha(\alpha (a)) && 
\\ & = \alpha(\alpha(a)+\alpha(z)) && (\text{Eq. \eqref{un}})
\\ & = \alpha(\alpha(a)\cdot z)\cdot z && (\text{Lemma \ref{lemma: Lukas near semi sono integrali}}) 
\\ & = \alpha(z)\cdot z && (\text{Eq. \eqref{du}}) 
\\ & = 0
\end{align*}
Therefore $ a = (x \cdot \alpha(y)) + (y\cdot \alpha(x)) = 0 $.  Since $ \langle R, + \rangle $ is a join-semilattice with 0 as least element, $ (x \cdot \alpha(y)) + (y\cdot \alpha(x)) = 0 $ implies that   $x \cdot \alpha(y) = 0 $ and $ y\cdot \alpha(x) = 0 $. Using Lemma \ref{lem: x leq y iff x * y' = 0}, we get $ x\leq y $ and $ y\leq x $, proving that $ x = y $ as desired.
\end{proof}
\end{theorem}

\begin{theorem}\label{th: near Luk semirings are arithmetical}
The variety of \L ukasiewicz near semirings is arithmetical, with witness Mal'cev term
$$ p(x,y,z) = \alpha ((\alpha(x\cdot\alpha (y))\cdot\alpha(z)) + (\alpha(z\cdot\alpha (y))\cdot\alpha(x))). $$
\begin{proof}
We first show that the term $ p(x,y,z)$ is a Mal'cev term for the variety of \L ukasiewicz near semiring: $ p(x,y,y) = x $ and $ p(x,x,y) = y $. 
\begin{align*}
& p(x,y,y) = \alpha ((\alpha(x\cdot\alpha (y))\cdot\alpha(y)) +(\alpha(y\cdot\alpha (y))\cdot\alpha(x)))
\\ & = \alpha ((\alpha(x + y) + \alpha (x)) & (\text{Lemma \ref{lemma: Lukas near semi sono integrali}})
\\ & = \alpha (\alpha (x)) = x  & (\text{Lemma \ref{lemma: aritmetica}})
\end{align*}
Similarly, 
\begin{align*}
& p(x,x,y) = \alpha ((\alpha(x\cdot\alpha (x))\cdot\alpha(y)) +(\alpha(y\cdot\alpha (x))\cdot\alpha(x)))
\\ & = \alpha ((\alpha(y) + \alpha (x+y)) & (\text{Lemma \ref{lemma: Lukas near semi sono integrali}})
\\ & = \alpha (\alpha (y)) = y  & (\text{Lemma \ref{lemma: aritmetica}})
\end{align*}
Therefore the variety of \L ukasiewicz near semirings is congruence permutable. Moreover, the following ternary term
$$ M(x,y,z) =  \alpha(\alpha (x) + \alpha (y)) + \alpha(\alpha (y) + \alpha (z)) + \alpha(\alpha(z)+\alpha(x)) $$
is a majority term for the variety of \L ukasiewicz near semiring. A simple calculation shows that $ M(x, x, y) = M(x, y, x) = M(y, x, x) = x
 $. This proves that the variety considered is also congruence distributive, hence by definition it is arithmetical as claimed. 
\end{proof}
\end{theorem}

%
%

\section{Orthomodular lattices as near semirings}\label{sec: 5}

Orthomodular lattices were introduced in 1936 by Birkhoff and von Neumann as an algebraic account of the logic of quantum mechanics. A detailed discussion can be found in \cite{Beran11, Kalmbach83}. The aim of this section is to provide a representation of orthomodular lattices as involutive near semirings.

Let us briefly recall that an \textit{orthomodular lattice} (\emph{OML}, for short) is an algebra $ \mathbf{L}=\langle L, \vee, \wedge, ', 0, 1\rangle $ of type $ \langle 2,2,1,0,0\rangle $ such that $ \langle L, \vee, \wedge, 0, 1\rangle  $ is a bounded lattice, $ ' $ is an orthocomplementation, i.e. $ x\wedge x' = 0 $, $ x\vee x' = 1 $. Furthermore $ ' $ is an involutive, antitone map ($ x\leq y $ implies $ y'\leq x' $) that satisfies the so called \textit{orthomodular law}:
\begin{equation}\label{eq: orthomodular law}
 x\leq y \;\;\Rightarrow\;\;  y= x\vee (y\wedge x').
\end{equation}
\noindent
The orthomodular law can be equivalently expressed by the identity
\begin{equation}\label{eq: orthomodular law (identity)}
(x\vee y)\wedge (x\vee (x\vee y)') = x,
\end{equation}
\noindent
which, in turn, is equivalent to the dual form:
\begin{equation}\label{eq: orthomodular law (dual identity)}
(x\wedge y)\vee (y\wedge (x\wedge y)') = y.
\end{equation}
In the next lemma we recap some basic facts relative to OMLs which will be useful in what follows. Let $ a,b $ two elements of an OML $ \mathbf{L} $, we say that $ a $ and $ b $ \textit{commute} (in symbols $ a C b $) iff $ a=(a\wedge b)\vee (a\wedge b') $. For the proof of Lemma \ref{lem: facts about OML}, see \cite{Beran11} or \cite{Kalmbach83}.
\begin{lemma}\label{lem: facts about OML}
Let $ \mathbf{L} $ an orthomodular lattice and $ a,b,c\in L $. Then 
\begin{itemize}
\item[(i)] If $ aCb $ then $ bCa $
\item[(ii)] If $ a\leq b $ then $ aCb $
\item[(iii)] If $ aCb $ then $ aCb' $
\item[(iv)] If two elements among $ a,b,c $ commutes with the third, then $ (a\vee b)\wedge c = (a\wedge c)\vee (a\wedge c) $ and $ (a\wedge b)\vee c = (a\vee c)\wedge (b\vee c) $
\end{itemize}
\end{lemma}
\noindent
In the previous section, we introduced  \L ukasiewicz near semirings to represent basic algebras. Here, to provide a similar representation of OMLs, we will consider orthomodular near semirings.
\begin{definition}\label{def: orthomodular near semirings}
An \emph{orthomodular near semiring} $ \mathbf{R} $ is a \L ukasiewicz near semiring that fulfills the following identity:
\begin{equation}\label{eq: def orthomodular near semirings}
x = x\cdot (x + y) 
\end{equation}
\end{definition}

The next lemma shows some basic properties of orthomodular near semirings.

\begin{lemma}\label{lem: aritmetica ONS}
Let $ \mathbf{R} $ be an orthomodular near semiring. Then: 
\begin{itemize}
\item[(a)] $ x \cdot x = x $;
\item[(b)] $ x = x\cdot \alpha((\alpha(y \cdot \alpha(x)) \cdot \alpha(x)) $;
\item[(c)] $ x + \alpha(x) = 1 $;
\item[(d)] If $ x\leq y $ then $ x\cdot y = y $.
\end{itemize}
\begin{proof}
(a) Straightforward, by setting $ y = 0 $ (or also $ x = y $) in equation \eqref{eq: def orthomodular near semirings}.  \\
(b) follows directly using equation \eqref{eq: def orthomodular near semirings} and Lemma \ref{lemma: Lukas near semi sono integrali}-(e). \\
(c) By Lemma \ref{lemma: Lukas near semi sono integrali}-(e), we have $ x + \alpha(x) = \alpha(\alpha(x\cdot x)\cdot x) = \alpha (\alpha (x)\cdot x) = 1 $, where we have used (a). \\
(d) Let $ a\leq b $, then $ a + b = b $. Therefore $ a = a\cdot (a + b) = a\cdot b $.     
\end{proof} 
\end{lemma}

We first show that an orthomodular near semiring can always be obtained out of an OML. 
\begin{theorem}\label{th: da OML a ONS}
Let $\mathbf{L}=\langle L, \vee , \wedge , ^{'}, 0, 1\rangle $ an orthomodular lattice and define multiplication via the so-called {\it Sasaki projection}: $ x\cdot y:=(x\vee y')\wedge y $. Then $ \mathbf{R}(\mathbf{L})= \langle L, + , \cdot , ^{'} , 0,1\rangle $ is an orthomodular near semiring,  where $x+y=x\lor y$.
\begin{proof}
It is evident that $ \langle L, \vee, 0\rangle $ is a commutative, idempotent monoid. Furthermore, $ x\cdot 1 = (x\vee 1')\wedge 1 = (x\vee 0)\wedge 1 = x $, and $ 1\cdot x = (1\vee x')\wedge x = 1\wedge x = x $. Therefore $ \langle R, \cdot , 1\rangle $ is a groupoid with $ 1 $ as neutral element. To prove right distributivity we make use of Lemma \ref{lem: facts about OML}. Upon observing that $ z'\leq x\vee z' $, $ z'\leq y\vee z' $, we have that $ z' $ {commutes} (in the sense of Lemma \ref{lem: facts about OML}) with both $ x\vee z' $ and $ y\vee y' $, therefore $ z $ does. For this reason we get: 
\begin{align*}
& (x\vee y)\cdot z = ((x\vee y)\vee z')\wedge z && (\text{Definition}) 
\\ & = ((x\vee z')\vee(y\vee z')) \wedge z && (\text{Lattice prop.})
\\ & = ((x\vee z')\wedge z)\vee ((y\vee z')\wedge z) && (\text{Lemma \ref{lem: facts about OML}-(iv)}) 
\\ & = (x\cdot z)\vee (y\cdot z). 
\end{align*}
It is not difficult to check that $ 0 $ annihilates multiplication. Indeed, $ x\cdot 0 = (x\vee 0')\wedge 0 = 0 $ and $ 0\cdot x = (0\vee x')\wedge x = x'\wedge x = 0 $. We now show that $ \mathbf{R}(\mathbf{L}) $ is \L ukasiewicz near semiring (see Definition \ref{def: Luka near semiring}). First let us observe that: 
\begin{align*}
& (x\cdot y')' \cdot y' = (((x\vee y)\wedge y')'\vee y) \wedge y' && (\text{Definition})
\\ & = (((x\vee y)' \vee y)\vee y) \wedge y' && (\text{De Morgan}) 
\\ & = ((x\vee y)' \vee y)\wedge y' && (\text{Ass., Idem.})
\end{align*}
Reasoning similarly one gets $  (y\cdot x')' \cdot x' = ((x\vee y)' \vee x)\wedge x' $. \\
Simply observing that $ x\leq x\vee y $ and applying the orthomodular law, we have $ x\vee y= x\vee ((x\vee y)\wedge x') $. Therefore,
\begin{align*}
& (x\vee y)'= (x\vee ((x\vee y)\wedge x'))'  &&
\\ & = x' \wedge ((x\vee y)\wedge x')' && (\text{De Morgan})   
\\ & = x' \wedge ((x\vee y)'\vee x) && (\text{De Morgan})   
\\ & = ((x\vee y)' \vee x)\wedge x' && (\text{Comm.}) 
\\ & = (y\cdot x')' \cdot x'
\end{align*}
Analogously, using the fact that $ y\leq x\vee y $ one gets, by the orthomodular law, that $ (x\vee y)' = ((x\vee y)' \vee y)\wedge y' = (x\cdot y')' \cdot y' $. Therefore $ (x\cdot y')' \cdot y' = (y\cdot x')' \cdot x' $ as claimed. We finally check that also equation \eqref{eq: def orthomodular near semirings} holds. This is a simple consequence of the orthomodular law: $ x\cdot (x+y) = (x\vee (x\vee y)')\wedge (x\vee y) = x $ by equation \eqref{eq: orthomodular law (identity)}. Therefore, $ \mathbf{R}(\mathbf{L})= \langle L, \vee , \cdot , ^{'} , 0,1\rangle $ is an orthomodular near semiring.\footnote{Notice that $ \langle L, \vee ,\cdot\rangle $, in general, is not a lattice.}
\end{proof}
\end{theorem}
 \noindent
 We can also prove the converse, stating a correspondence between orthomodular lattices and orthomodular near semirings. 
 
\begin{theorem}\label{th: da orthomodular near semirings a orthmodular lattices}
Let $ \mathbf{R} $ be an orthomodular near semiring. Setting $ x\vee y = x + y $, $ x' = \alpha (x) $, and then defining $ x\wedge y = (x' \vee y')' $, then $ \mathbf{L}(\mathbf{R})= \langle R, \vee , \wedge , ', 0,1\rangle $ is an orthomodular lattice. 
\begin{proof}
Since $ \mathbf{R} $ is integral we know that $ \langle R, +\rangle $ is a join-semilattice with 1 as top element, and consequently $ \langle R, \vee\rangle $ is. On the other hand, since $ \alpha $ is an antitone involution then $ \langle R, \wedge\rangle $ is a meet-semilattice with 0 as least element. As meet and join are defined dually, $ \langle R, \vee, \wedge, 0, 1\rangle $ is a bounded lattice. Furthermore, $ x\vee x' = 1 $ is guaranteed by Lemma \ref{lem: aritmetica ONS} and thus it follows that $ x\wedge x' = 0 $. \\
We are left with the task of showing that the orthomodular law holds too. So, suppose $ a\leq b $, then 
\begin{align*}
& a \vee (b\wedge a') = a + \alpha(\alpha (b) + a) && 
\\ & = a + \alpha(a + \alpha (b)) && (\text{Comm.}) 
\\ & = a + (\alpha(a\cdot b) \cdot b) && (\text{Lemma \ref{lemma: Lukas near semi sono integrali}})
\\ & = (a\cdot b) + (\alpha(a)\cdot b) && (\text{Lemma \ref{lem: aritmetica ONS}-(d)})
\\ & = (a + \alpha(a)) \cdot b && (\text{Distr.})
\\ & = 1 \cdot b = b.
\end{align*}
This allows to conclude that $ \mathbf{L}(\mathbf{R}) $ is an orthomodular lattice. 
\end{proof}
\end{theorem}
\noindent
The theorems above have shown how to get an orthomodular lattice out of an orthomodular semirings and viceversa. In other words, there are maps $ f ,\, g $, from the variety of orthomodular lattices to the variety of orthomodular semirings and from the variety of orthomodular near semirings to the variety of orthomodular lattices, respectively, assigning to any OML an orthomodular semiring, and vice versa. We now show that:
\begin{theorem}\label{th: OML mappe mutually inverse}
The maps $ f $ and $ g $ are mutually inverse: $ \mathbf{L}=\mathbf{L}(\mathbf{R}(\mathbf{L})) $ and $ \mathbf{R}=\mathbf{R}(\mathbf{L}(\mathbf{R})) $. 
\begin{proof}
Let $ \mathbf{L} \langle L, {\vee}, {\wedge}, ', 0, 1\rangle $ be an orthomodular lattice. It follows from Theorem \ref{th: da OML a ONS} that $ \mathbf{R}(\mathbf{L}) $ is an orthomodular near semiring, and from Theorem \ref{th: da orthomodular near semirings a orthmodular lattices} that the structure $\mathbf{L}(\mathbf{R}(\mathbf{L})) =\langle L, \bar{\vee}, \bar{\wedge}, ', 0, 1\rangle $ is an orthomodular lattice. It is straightforward to check that the involutions on $ \mathbf{L}(\mathbf{R}(\mathbf{L})) $ and $ \mathbf{L} $ coincide, as well as $ x\bar{\vee} y = x\vee y $. Therefore we also have that $ x\bar{\wedge} y = (x' \bar{\vee} y')' = (x' \vee y')' = x\wedge y $. So  $ \mathbf{L}=\mathbf{L}(\mathbf{R}(\mathbf{L})) $. \\
On the other hand, by Theorems  \ref{th: da OML a ONS} and \ref{th: da orthomodular near semirings a orthmodular lattices} we obtain that the structure $ \mathbf{R}(\mathbf{L}(\mathbf{R}))  =\langle R, \widehat{+}, \widehat{\cdot}, ^{\widehat{\alpha}}, 0, 1\rangle $ is an orthomodular near semiring.  Again it is straightforward to check that $ \widehat{\alpha} (x) = \alpha (x) $ and $ x\;\widehat{+}\; y = x + y $. It is less evident that $ x\;\widehat{\cdot} \; y = x\cdot y $. Indeed: 
\begin{align*}
& x\;\widehat{\cdot}\; y = (x\vee y')\wedge y && 
\\ & = ((x\vee y')'\vee y')' &&
\\ & = \alpha(\alpha(x + \alpha(y)) +\alpha( y)),  &&
\end{align*}
where $ \vee $, $ \wedge $ and $ ^{'} $ are join, meet and complementation, respectively, of the orthomodular lattice $ \mathbf{L}(\mathbf{R}) $. We are finally left with showing that $\alpha(\alpha(x +\alpha( y)) +\alpha( y)) = x\cdot y $. \\
\noindent Our first move is to prove that $\alpha(\alpha(x +\alpha( y)) +\alpha( y)) = y\cdot (x + \alpha(y))$.
\begin{align*}
& \alpha(\alpha(x +\alpha( y)) +\alpha( y))  = \alpha(\alpha( y)+\alpha(x +\alpha( y)))  
\\ & = \alpha(\alpha(y)\cdot(x+\alpha(y)))\cdot(x+\alpha(y)) && (\text{Lemma \ref{lemma: Lukas near semi sono integrali}})
\\ & = \alpha(\alpha (y))\cdot(x+\alpha(y)) && \eqref{eq: def orthomodular near semirings}
\\ & = y\cdot(x+\alpha(y)).
\end{align*}
Moreover,
\begin{align*}
& \alpha((x \cdot y)\cdot y)\cdot y =  \alpha((x \cdot y)+\alpha(y)) && (\text{Lemma \ref{lemma: Lukas near semi sono integrali}})
\\ & =   \alpha(\alpha(y)+(x \cdot y)) && (\text{Comm.}) 
\\ & = \alpha(\alpha(y)\cdot \alpha(x \cdot y))\cdot \alpha(x\cdot y) && (\text{Lemma \ref{lemma: Lukas near semi sono integrali}}) 
\\ & = \alpha(\alpha(y))\cdot \alpha(x\cdot y) && (\text{by Lemma \ref{lem: aritmetica ONS}, since $\alpha(y)\leq \alpha(x \cdot y)$}) 
\\ & = y\cdot \alpha(x\cdot y).
\end{align*}
Using the derivation above, which we will refer to as ($ \star $), we finally prove our claim:
\begin{align*}
& y\cdot (x + \alpha(y)) = y\cdot \alpha(\alpha(x \cdot y)\cdot y) && (\text{Lemma \ref{lemma: Lukas near semi sono integrali}}) 
\\ & = \alpha((\alpha(x\cdot y)\cdot y)\cdot y)\cdot y && (\text{$ \star $})
\\ & = \alpha(\alpha(x\cdot y)\cdot y)\cdot y && (\text{by Lemma \ref{lem: aritmetica ONS}, since $\alpha(x\cdot y)\cdot y\leq y$}) 
\\ & = (x+\alpha(y))\cdot y && (\text{Lemma \ref{lemma: Lukas near semi sono integrali}})
\\ & = (x\cdot y)+(\alpha(y)\cdot y) && (\text{Right Distr.})
\\ & = x\cdot y. && (\text{Lemma \ref{lemma: Lukas near semi sono integrali}})
\end{align*}

\end{proof}
\end{theorem}


\section{Central elements and decomposition}\label{sec: 6}

The aim of this section is to give a a characterization of the central elements and consequently some decomposition theorems for the variety of integral involutive near semirings. Such results apply to both the variety of \L ukasiewicz near semirings and orthomodular near semirings as the they are both integral. The section relies on the ideas developed in \cite{Sal} and \cite{Ledda13} on the general theory of \textit{Church algebras}.\footnote{Diverse applications of this theory can be found in \cite{CGKGLP, Bonzio15}.}

The notion of Church algebra is based on the simple observation that many well-known algebras, including Heyting algebras, rings with unit and combinatory algebras, possess a ternary term operation $ q $, satisfying the equations: $ q(1,x,y)= x $ and $ q(0,x,y)= y $. The term operation $ q $ simulates the behavior of the if-then-else connective and, surprisingly enough, this yields rather strong algebraic properties. 

An algebra $\mathbf{A}$\ of type $\nu$\ is a \emph{Church algebra}\ if there are term definable elements $0^{\mathbf{A}},1^{\mathbf{A}}\in A$\ and a ternary term operation $q^{\mathbf{A}}$\ s.t., for all $a,b\in A$, $q^{\mathbf{A}}\left(1^{\mathbf{A}},a,b\right)  =a$\ and $q^{\mathbf{A}}\left(  0^{\mathbf{A}},a,b\right)  =b$. A variety $\mathcal{V}$\ of type $\nu$\ is a Church
variety\ if every member of $\mathcal{V}$\ is a Church algebra with respect to the same term $q\left(  x,y,z\right)  $\ and the same constants $0,1$. 

Taking up a suggestion from Diego Vaggione \cite{Vaggio}, we say that an element \textit{e} of a Church algebra \textbf{A} is \textit{central} if the pair $ (\theta(e,0),\theta (e,1) ) $ is a pair of factor congruences on \textbf{A}. A central element $ e $ is nontrivial when $ e\not\in\{0,1\} $. We denote the set of central elements of \textbf{A} (the {\it centre}) by $\mathrm {Ce}({A})$.

Setting
$$ x\wedge y=q(x,y,0),\;\; x\vee y= q(x,1,y)\;\; x^{*}=q(x,0,1) $$
we can report the following general result on Church algebras:
\begin{theorem}\label{th: Boolean algebra of centrals}
\emph{\cite{Sal}} Let $ \mathbf{A} $ be a Church algebra. Then 
\[
\mathrm {Ce}(\mathbf{A})=\langle \mathrm {Ce}(A),\wedge,\vee,^{*},0,1\rangle
\]
is a Boolean algebra which is isomorphic to the Boolean algebra of factor congruences of $\mathbf{A}$.
\end{theorem}


If $\mathbf{A}$ is a Church algebra of type $\nu$ and $e\in A$ is a central
element, then we define $\mathbf{A}_{e}=(A_{e},g_{e})_{g\in\nu}$ to be the
$\nu$-algebra defined as follows:

\begin{equation}\label{eq:opAe}
A_{e}=\{e\wedge b:b\in A\};\quad g_{e}(e\wedge\overline{b})=e\wedge
g(e\wedge\overline{b}), 
\end{equation}
where $ \overline{b} $ denotes the a n-tuple $ b_1,...,b_n $ and $e\wedge\overline{b} $ is an abbreviation for $ e\wedge b_1,...,e\wedge b_n $.

By \cite[Theorem 4]{Ledda13}, we have that:

\begin{theorem}\label{th: decomposizione Church algebras}
\label{relat}Let $\mathbf{A}$ be a Church algebra of type $\nu$ and $e$ be a
central element. Then we have:

\begin{enumerate}
\item For every $n$-ary $g\in\nu$ and every sequence of elements $\overline
{b}\in A^{n}$, $e\wedge g(\overline{b})=e\wedge g(e\wedge\overline{b})$, so
that the function $h:A\rightarrow A_{e}$, defined by $h(b)=e\wedge b $, is a
homomorphism from $\mathbf{A}$ onto $\mathbf{A}_{e}$.

\item $\mathbf{A}_{e}$ is isomorphic to $\mathbf{A}/\theta(e,1)$. It follows
that $\mathbf{A}=\mathbf{A}_{e}\times\mathbf{A}_{e^{\prime}}$ for every
central element $e$, as in the Boolean case.
\end{enumerate}
\end{theorem}
\noindent

\begin{proposition}\label{prop: integral near semiring are a Church variety}
The class of intergral involutive near semirings is a Church variety, as witnessed by the term: 
\[
q(x,y,z) = (x\cdot y)+ (\alpha(x)\cdot z).
\]
\begin{proof}
Suppose $ \mathbf{R} $ is an integral involutive near semiring and $ a,b\in R $. Then $ q(1,a,b)=(1\cdot a)+(\alpha(1)\cdot b) = a + (0\cdot b) = a+0 = a. $ and $ q(0,a,b)=(0\cdot a)+(\alpha(0)\cdot b) = 0 + (1\cdot b) = 0+b = b  $. 
\end{proof}
\end{proposition}
\noindent
Since both the varieties of \L ukasiewicz and orthomodular near semirings are subvarieties of integral involutive near semiring, it follows that both of them are Church varieties. In this section we apply the theory of Church algebras to the more general class of integral involutive near semirings.  
According with the results in \cite[Proposition 3.6]{Sal}, in a Church variety central elements are amenable to a very general description.
\begin{proposition}\label{prop: description of central elements in Church varieties}
If $\mb{A}$ is a Church algebra of type $ \nu $ and $ e\in A $, the following conditions are equivalent: 
\begin{itemize}
\item[(1)] e is central;
\item[(2)] for all $ a,b, \vec{a},\vec{b}\in A $: 
\begin{itemize}
\item[\textbf{a)}] $ q(e,a,a)=a $,
\item[\textbf{b)}] $ q(e,q(e,a,b),c)=q(e,a,c)=q(e,a,q(e,b,c)) $,
\item[\textbf{c)}] $ q(e,f(\vec{a}),f(\vec{b}))= f(q(e,a_{1},b_{1}),...,q(e,a_{n},b_{n})) $, for every $ f\in\nu $, 
\item[\textbf{d)}] $ q(e,1,0)=e $.
\end{itemize}
\end{itemize}
\end{proposition}
\noindent
In case \textbf{A} is an integral involutive near semiring, condition (a) reduces to 
\begin{equation}\label{eq: centrali (a)}
(e\cdot a)+ (\alpha(e)\cdot a) = a.
\end{equation}
Conditions (b) read
\begin{equation}\label{eq: centrali (b1)}
(e\cdot ((e\cdot a) + (\alpha(e)\cdot b)))+(\alpha(e)\cdot c))= (e\cdot a)+(\alpha(e) \cdot c),
\end{equation}
\begin{equation}\label{eq: centrali (b2)} 
 (e\cdot a)+(\alpha(e) \cdot c) = (e\cdot a) + (\alpha(e) \cdot ((e\cdot b) + (\alpha(e)\cdot c))). 
\end{equation}
\noindent
Condition (c), whenever $ f $ is the constant 0, expresses a property that holds for every element: $ (e\cdot 0) + (\alpha(e) \cdot 0) = 0 $. On the other hand, if $ f $ coincides with the nullary operation $ 1 $, we obtain (for a central element $e$)
\begin{equation}\label{eq: e+e' = 1 sui centrali}
q(e,1,1)= (e \cdot 1) + (\alpha(e)\cdot1)= e + \alpha(e) = 1.
\end{equation}
If $ f $ coincides with the involution, (c) reads
\begin{equation}\label{eq: centrali (c1)}
(e\cdot \alpha(a))+ (\alpha(e)\cdot \alpha(b))=\alpha((e\cdot a) + (\alpha(e)\cdot b)). 
\end{equation}
\noindent
Whenever $ f $ is $+$, we obtain: 
\begin{equation}\label{eq: centrali (c2)}
(e\cdot(a+c)) + (\alpha(e)\cdot(b+d))= ((e\cdot a)+ (\alpha(e)\cdot b))+((e\cdot c) + (\alpha(e)\cdot d)), 
\end{equation}
this, by the associativity of the sum, is equal to
\begin{equation}\label{eq: sorta di associatività per i centrali}
((e\cdot a) + (e\cdot c)) + ((\alpha(e)\cdot b) +  (\alpha(e)\cdot d)),
\end{equation}
which is a sort of distributivity restricted to central elements. Whenever $ f $ is the multiplication, this condition simplifies to
\begin{equation}\label{eq: centrali (c3)}
(e\cdot(a\cdot c)) + (\alpha(e)\cdot(b\cdot d))= ((e\cdot a)+ (\alpha(e)\cdot b))\cdot ((e\cdot c) + (\alpha(e)\cdot d)).
\end{equation}
\noindent
Condition (d) expresses a general property that holds true for every element: $ (e\cdot 1) + (\alpha (e) \cdot 0) = e + 0 = e $. 
We have just seen in Proposition \ref{prop: description of central elements in Church varieties} that, in Church algebras, central elements can be described by means of identities. This, in fact, will be very useful in proving the results in this section. However, we aim to show that the axiomatization of central elements can be streamlined to a \emph{minimal set} (see Appendix \ref{sec:app}) of two identities only. The next lemma introduces some results which are very useful to prove the minimality of such an axiomatization.  
\begin{lemma}\label{lemma: utile per i centrali}
Let $ \mathbf{R} $ be an integral involutive near semiring, and $ e\in R $ an element that satisfies the following identities: 
\begin{enumerate}
\item $ (e\cdot \alpha(x))+ (\alpha(e)\cdot \alpha(y))=\alpha((e\cdot x) + (\alpha(e)\cdot y)) $;
\item $ (e\cdot(x\cdot z)) + (\alpha(e)\cdot(y\cdot u))= ((e\cdot x)+ (\alpha(e)\cdot y))\cdot ((e\cdot z) + (\alpha(e)\cdot u)) $.
\end{enumerate}
 Then $e$ satisfies the following:
\begin{itemize}
\item[(i)] $ (e\cdot x) + \alpha (e) = x + \alpha (e) $;
\item[(ii)] $ e\cdot (e\cdot x) = e\cdot x = (e\cdot x)\cdot e $;
\item[(iii)] $ e\cdot\alpha (e) = 0 $;
\item[(iv)] $ e\cdot x = x\cdot e $; 
\item[(v)] $ e\cdot (x + y) = (e\cdot x) + (e\cdot y) $;
\item[(vi)] if $ x\leq y $ then $ e\cdot x\leq e\cdot y $;
\item[(vii)] $ e\cdot (\alpha(e)\cdot x) = 0 $.
\end{itemize}
\begin{proof}
(i) Since $ e\leq 1 $, then $ e\cdot x\leq 1\cdot x = x $. Therefore $ (e\cdot x) + \alpha (e) \leq x+\alpha (e) $. For the converse, first notice that, as $ e\cdot \alpha(x)\leq\alpha(x) $, then $ x\leq\alpha (e\cdot \alpha(x)) = (e\cdot x) + \alpha (e) $, where the last equality is obtained by setting $ y = 1 $ in identity (1) (and the fact that $ \alpha $ is an involution). \\
(ii) The first equality readily follows from (2) upon setting $ y = u = 0 $ and $ x = 1$, while the second by setting $ y = u = 0 $ and $ z = 1 $. \\
(iii) can be derived by setting $ x = u = 1 $ and $ y = z = 0 $ in identity (2). \\
(iv) 
\begin{align*}
& e\cdot x = (e\cdot x)\cdot e && (\text{ii})
\\ & = ((e\cdot x)\cdot e) + (\alpha(e)\cdot e) && (\text{iii})
\\ & = ((e\cdot x) + \alpha(e))\cdot e && (\text{Distr})
\\ & = (x+\alpha (e))\cdot e && (\text{i}) 
\\ & = (x\cdot e) + (\alpha (e)\cdot e) && (\text{Distr})
\\ & = (x\cdot e) + 0 && (\text{iii}) 
\\ & = x\cdot e. 
\end{align*}
(v) 
\begin{align*}
& e\cdot (x + y) = (x+y)\cdot e && (\text{iv})
\\ & = (x\cdot e) + (y\cdot e) && (\text{Distr})
\\ & = (e\cdot x) + (e\cdot y) && (\text{iv})
\end{align*}  
(vi) Let $ x\leq y $, i.e. $ x + y = y $. Then $ e\cdot y = e\cdot (x +y) = (e\cdot x) + (e\cdot y) $, i.e. $ e\cdot x \leq e\cdot y $.  \\
(vii) In case  $ y = u = 0 $, in condition (3), we obtain: $e\cdot(x\cdot z)=(e\cdot x) \cdot (x \cdot y) $. 
  If, moreover, $ x = \alpha(e) $, we obtain that $ e\cdot (\alpha(e)\cdot z) = (e\cdot \alpha(e))\cdot (e\cdot z) = 0 $, by (iii).
\end{proof}
\end{lemma}
\noindent
We now put Lemma \ref{lemma: utile per i centrali} to good use and prove that, in an involutive near semiring, central elements are neatly characterized by two simple equations.

\begin{theorem}\label{th: centrali nei near semirings}
Let $ \mathbf{R} $ be an involutive near semiring. Then an element $ e\in R $ is central if and only if it satisfies the following equations for any $ x,y,z,u\in R $:
\begin{enumerate}
\item $ (e\cdot \alpha(x))+ (\alpha(e)\cdot \alpha(y))=\alpha((e\cdot x) + (\alpha(e)\cdot y)) $;
\item $ (e\cdot(x\cdot z)) + (\alpha(e)\cdot(y\cdot u))= ((e\cdot x)+ (\alpha(e)\cdot y))\cdot ((e\cdot z) + (\alpha(e)\cdot u)) $.
\end{enumerate}
\begin{proof}
($\Rightarrow $) If $ e $ is a central element then (1), (2) hold by Proposition \ref{prop: description of central elements in Church varieties}. \\
($\Leftarrow $) Using again Proposition \ref{prop: description of central elements in Church varieties}, and identities (1) and (2), we have to derive equations \eqref{eq: centrali (a)}, \eqref{eq: centrali (b1)}, \eqref{eq: centrali (b2)}, \eqref{eq: e+e' = 1 sui centrali} and \eqref{eq: centrali (c2)}. We start by deriving \eqref{eq: e+e' = 1 sui centrali}: upon setting $ x = y = 0$, identity (1) reads: $ e+\alpha(e) = \alpha (0) = 1 $.
Using \eqref{eq: e+e' = 1 sui centrali}, we obtain \eqref{eq: centrali (a)} as follows
\begin{align*}
& (e\cdot x)+ (\alpha(e)\cdot x)= (e + \alpha(e)) \cdot x && (\text{Distr.}) 
\\ & = 1\cdot x && \eqref{eq: e+e' = 1 sui centrali}
\\ & = x. 
\end{align*}
Equation \eqref{eq: centrali (c2)} immediately follows from the associativity of the sum and the fact that $ e\cdot (x + y) = (e\cdot x) + (e\cdot y) $ from Lemma \ref{lemma: utile per i centrali}.
In order to prove \eqref{eq: centrali (b1)} and \eqref{eq: centrali (b2)} we use some auxiliary facts stated in Lemma \ref{lemma: utile per i centrali}. 
\begin{align*}
& (e\cdot ((e\cdot a) + (\alpha(e)\cdot b))) + (\alpha(e)\cdot c) =
\\ & = (e\cdot (e\cdot a)) + (e\cdot (\alpha(e)\cdot b)) + (\alpha(e)\cdot c)  && (\text{Lemma \ref{lemma: utile per i centrali}.(v)})
\\ & = (e\cdot a) + (e\cdot (\alpha(e)\cdot b)) + (\alpha(e)\cdot c)  && (\text{Lemma \ref{lemma: utile per i centrali}.(ii)}) 
\\ & = (e\cdot a) + 0 + (\alpha(e)\cdot c)  && (\text{Lemma \ref{lemma: utile per i centrali}.(vii)}) 
\\ & = (e\cdot a) + (\alpha(e)\cdot c) 
\end{align*}
With a slight modification of the reasoning above one can derive condition \eqref{eq: centrali (b2)}. 
\end{proof}
\end{theorem}
\noindent
The next proposition yields a more informative version of the general result stated in Theorem \ref{th: Boolean algebra of centrals}. 
\begin{proposition}\label{prop: Boolean algebra of centrals}
Let $\mb{R}$ be an integral involutive near semiring and $ \mathrm{Ce}({R})$ the set of central elements of $\mathbf{R}$. Then $ \mathrm{Ce}(\mathbf{R})=\langle \emph{Ce}(R), +, \cdot, \alpha, 0,1\rangle $ is a Boolean algebra.
\begin{proof}
By Theorem \ref{th: Boolean algebra of centrals}, $\mathrm {Ce}(\mathbf{R})=\langle \mathrm {Ce}(R),\wedge,\vee,^{*},0,1\rangle $ is a Boolean algebra, where $ \wedge,\vee,^{*} $ are defined as follows
$$ x\wedge y=q(x,y,0),\;\; x\vee y= q(x,1,y)\;\; x^{*}=q(x,0,1) $$
Using this result, we just check that, for central elements, $ \wedge,\vee,^{*} $ coincide with $ \cdot, + , \alpha $, respectively. We can easily obtain that $ x\wedge y = q(x,y,0) = (x\cdot y ) + (\alpha(x)\cdot 0) = x\cdot y $, and $ x^* = q(x,0,1) = (x \cdot 0) + (\alpha (x)\cdot 1) = \alpha (x ) $. \\
It only remains to show that $ x + y = \alpha (\alpha (x)\cdot \alpha (y)) $. Notice first that, by equation \eqref{eq: centrali (c1)}, with $ a= 0 $, $ b = y' $ and $ e = x $ (this is legitimated by the fact that we are only concerned with central elements), we have
\begin{align*}
& x + (\alpha (x) \cdot y) = \alpha(\alpha(x) \cdot\alpha (y)) && (\dagger)
\end{align*}
Since, for central elements, multiplication coincides with the Boolean meet, we have that $ \alpha (x)\cdot \alpha (y) \leq \alpha (x) $ and $ \alpha (x)\cdot \alpha (y) \leq \alpha (y) $. As $\alpha$ is antitone, $ x \leq \alpha (\alpha (x)\cdot \alpha (y)) $ and $ y \leq \alpha (\alpha (x)\cdot \alpha (y)) $, which implies that $ x + y \leq \alpha (\alpha (x)\cdot \alpha (y)) + \alpha (\alpha (x)\cdot \alpha (y)) = \alpha (\alpha (x)\cdot \alpha (y)) $. For the converse, $ \alpha (x) \cdot y\leq y $, so $ x +(\alpha (x) \cdot y)\leq x+ y $, i.e. $ \alpha (\alpha (x)\cdot \alpha (y))\leq x + y $, by $ (\dagger) $. This proves that $ x + y = x \vee y $.  
\end{proof}
\end{proposition}
\noindent
From the previous proposition we have that if $ \mathbf{R} $ is an integral involutive near semiring and $ e $ is a central element, then $ \alpha (e) $ is also central. Our next step will be proving a decomposition theorem for involutive intergral near semiring. Let $ e $ be a central element of an integral involutive near semiring $ \mathbf{R} $, and set 
$$ [0,e]=\{x : x\leq e \} $$
A complementation can be naturally defined on $ [0,e] $ by setting $ x^{e}= e\cdot \alpha (x) $. Then, upon considering the algebra $ \mathbf{[0,e]} = \langle [0,e], +, \cdot , ^{e}, 0, e\rangle $, we can prove the following:

\begin{theorem}\label{th: decomposition in intervals}
Let $ \mathbf{R} $ an integral involutive near semiring and $e$ a central element of $ \mathbf{R} $. Then $ \mathbf{R}\cong \mathbf{[0,e]}\times \mathbf{[0,e']} $
\begin{proof}
As $ \mathbf{R} $ is a Church algebra, it satisfies Theorem \ref{th: decomposizione Church algebras}, hence all we have to prove reduce to the following: 
\begin{itemize}
\item[(1)] $ R_{e}= [0,e] $
\item[(2)] for $ x, y \leq e $, $ x + y = e \wedge (x + y) $,  $ x\cdot y = e\wedge (x\cdot y) $  \\
and $ x^{e} = e\wedge \alpha (x) $.
\end{itemize}
(1) Suppose $ x \in R_{e} $, i.e. $ x = e\wedge b $ for some $ b\in R $. By definition of $ \wedge $, $ e\wedge b = q(e,b,0) = (e\cdot b) + (\alpha (e)\cdot 0) = e\cdot b $. Now, as $ b\leq 1 $, by Lemma \ref{lemma: utile per i centrali} we have that  $ e\cdot b\leq e\cdot 1 = e $, i.e. $ x\in[0,e] $, proving $ R_{e} \subseteq [0,e] $.  For the converse, suppose $ x\in [0,e] $, i.e. $ x\leq e $. We want to find an element $ b\in R $ such that $ x = e\wedge b $. First notice that, under the assumption that $e$ is central, it follows by Theorem \ref{th: centrali nei near semirings} and Lemma \ref{lemma: utile per i centrali} that $ \alpha (e)\cdot x = 0 $, which we use to prove that
\begin{align*}
& 0 = \alpha(e)\cdot e && 
\\ & = \alpha (e)\cdot (e + x) && (\text{Assumption})
\\ & = (\alpha(e)\cdot e) + (\alpha(e) \cdot x) && (\text{Lemma \ref{lemma: utile per i centrali}})
\\ & = 0 + \alpha (e) \cdot x && (\text{Lemma \ref{lemma: utile per i centrali}})
\\ & =  \alpha (e) \cdot x. 
\end{align*}
We use the fact above to show that $ e\cdot x = x $. Since, by equation \eqref{eq: e+e' = 1 sui centrali}, $ 1 =e + \alpha (e)  $, we have that $ x =( e + \alpha (e))\cdot x = (e\cdot x) + (\alpha (e)\cdot x) = (e\cdot x) + 0 = e\cdot x   $. 
Remembering that $ e\wedge b = q(e,b,0) = (e\cdot b) + (\alpha (e)\cdot 0) = e\cdot b $ and setting $ b = x + \alpha (e) $ we get
\begin{align*}
& e\wedge b = e\cdot b && 
\\ & = e\cdot (x + \alpha (e)) && (\text{subs})
\\ & = (e \cdot x) + (e\cdot\alpha (e)) && (\text{Lemma \ref{lemma: utile per i centrali}}) 
\\ & = (e\cdot x) + 0 && (\text{Prop \ref{prop: Boolean algebra of centrals}}) 
\\ & = e\cdot x = x. 
\end{align*}
Therefore, $ x $ can be expressed as the meet of $ e $ with an element of $ R $, showing that $ [0,e]\subseteq R_{e} $. \\
(2) In this part of the proof we make use of the following facts
$$ x\wedge y = q(x,y,0) = x\cdot y \;\; \text{and} \;\; \text{if} \; x\leq e, \text{then } e\cdot x = x $$
Let $ x,y\leq e $. Then $  e\wedge (x + y) = e\cdot (x + y) = x+y $. Similarly, $ e\wedge (x\cdot y) = e\cdot (x\cdot y) = x\cdot y $. Finally $ x^{e} = e\wedge \alpha (e) = e\cdot \alpha (e) $
\end{proof}
\end{theorem}
 \noindent
Taking advantage from the fact that, in  a Church algebra, central elements are equationally characterizable (Proposition \ref{prop: description of central elements in Church varieties} and Theorem \ref{th: centrali nei near semirings}), we can prove the following:
\begin{proposition}\label{prop: atomi di A e Ae}
Let $ \mathbf{R} $ be a involutive integral near semiring, $ e\in\mathrm{Ce}(\mathbf{R}) $ and $ c\in R_e $. Then 
$$ c\in\mathrm{Ce}({R}) \Leftrightarrow c\in\mathrm{Ce}({R}_e) $$ 

\proof 
($\Rightarrow $) By Theorem \ref{th: centrali nei near semirings}, central elements are described by equations. Furthermore,  by Theorem \ref{th: decomposizione Church algebras}, $ h:\mathbf{R}\rightarrow\mathbf{R}_e $ is an onto homomorphism such that for every $ x\in R_e $, $ h(x)=x $. The fact that equations are preserved by homomorphisms yields the desired conclusion. \\
($\Leftarrow $) Let us observe that, since central elements are characterized by equations and equations are preserved by direct products, if $ c_1 $ and $ c_2 $ are central elements of two integral involutive near semirings $ \mathbf{R}_1 $ and $ \mathbf{R}_2 $, then $ (c_1,c_2)\in\mathrm{Ce}(\mathbf{R}_1 \times\mathbf{R}_2) $. Suppose $ c\in\mathrm{Ce}({R}_{e}) $, the image of $ c $ under the isomorphism of Theorem \ref{th: decomposizione Church algebras} is $ (c,0) $. On the other hand, $ 0 $ is always central element, therefore we have that $ (c,0) $ is a central element in $ \mathbf{R}_e \times \mathbf{R}_{e'} $, implying that $ c\in\mathrm{Ce}(\mathbf{R}) $, as $ \mathbf{R}\cong\mathbf{R}_e\times\mathbf{R}_{e'} $. 
\endproof 
\end{proposition} 

We have seen, in Proposition \ref{prop: Boolean algebra of centrals}, that $ \mathrm{Ce}(\mathbf{R})=\langle \emph{Ce}(R), +, \cdot, \alpha, 0,1\rangle $ is a Boolean algebra. Therefore it makes sense to consider the set of its atoms, which we denote by $At(\mathbf{R}) $.
\begin{lemma}\label{lemma: atomi}
If $ \mathbf{R} $ is an involutive integral near semiring and $ e\in At(\mathbf{R}) $, an atomic central element of $ \mathbf{R} $, then $ At(\mathbf{R}_{\alpha(e)})=At(\mathbf{R})\setminus\{e\} $.
\proof
($ \supseteq $) Suppose that $ e $ is an atom of the Boolean algebra $ \mathrm{Ce}(\mathbf{R}) $. Then, for any other atomic central element $ c\in\mathbf{R} $, $ c\wedge e = c\cdot e = e\cdot c = 0 $, therefore $ \alpha(e) + \alpha (c) =1 $. Furthermore, $ c = 1\cdot c = (e + \alpha (e))\cdot c = (e\cdot c) + (\alpha (e)\cdot c) = 0 + (\alpha(e)\cdot c) = \alpha (e)\cdot c $, which shows that $ c\leq \alpha (e) $. Thus, by Proposition \ref{prop: atomi di A e Ae}, $ c\in\mathbf{R}_{\alpha (e)} $. We have to show that $ c $ is also an atom. So, suppose $ d $ is a central element of $ \mathbf{R}_{\alpha(e)} $ such that $ d <\ c $, then, by Proposition \ref{prop: atomi di A e Ae}, $ d $ is a central element of $ \mathbf{R} $ and as, by assumption, $ c\in At(\mathbf{R}) $, then necessarily $ d=0 $, showing that $ c $ is also an atom in $ \mathbf{R}_{\alpha(e)} $. \\
$ (\subseteq ) $ Suppose $ c\in At(\mathbf{R}_{\alpha(e)}) $, then in particular $ c $ is a central element of $ \mathbf{R}_{\alpha(e)} $ and, by Proposition \ref{prop: atomi di A e Ae}, $ c\in\mathrm{Ce}(\mathbf{R}) $. Let $ d\in\mathrm{Ce}(\mathbf{R}) $, with $ c <\ d $, then we have $ d\leq \alpha(e) $ and therefore $ d\in\mathrm{Ce}(\mathbf{R}_{\alpha(e)}) $ by Proposition \ref{prop: atomi di A e Ae}. As, by assumption, $ c\in At(\mathbf{R}_{\alpha(e)})$ then $ d=0 $, which shows that $ c $ is an atomic central. We finally claim that $ c\neq e $. Indeed, suppose by contradiction that $ c= e $, then since $ c\leq \alpha(e) $ we have $ e\leq \alpha(e) $, i.e. $ e=e\cdot \alpha(e) = 0 $ which is a contradiction, as $ e $ is atomic central by hypothesis.
\endproof
\end{lemma}

Lemma \ref{lemma: atomi} will be useful in proving the following 
\begin{theorem}\label{th:decomposizione in prodotto di d.i.} 
Let $ \mathbf{R} $ be an involutive integral near semiring such that $ \mathrm{Ce}(\mathbf{R}) $ is an atomic Boolean algebra with countably many atoms, then
$$ \mathbf{R}=\prod_{e\in At(\mathbf{R})} \mathbf{R}_e $$
is a decomposition of $ \mathbf{R} $ as a product of directly indecomposable algebras. 
\proof
The claim is proved by induction on the number of elements of $ At(\mathbf{R}) $. If $ 1 $ is the only central atomic element, then $ \mathbf{R} $ is directly indecomposable and clearly $ \mathbf{R}=\mathbf{R}_1 $. If there is an atomic central element $ e\neq 1 $, then $ \mathbf{R}=\mathbf{R}_e\times\mathbf{R}_{\alpha(e)} $ by Theorem \ref{th: decomposizione Church algebras}. On the other hand $ \mathrm{Ce}(\mathbf{R}_{e})=\{0,e\} $, because if $ \mathbf{R}_e $ had another element, say $ d $, then $ d $ would be a central element of $ \mathbf{R} $ in virtue of Proposition \ref{prop: atomi di A e Ae} and $ 0 <\ d <\ e $ contradicting the fact that $ e $ is an atom. Consequently $ \mathbf{R}_e $ is directly indecomposable. By Lemma \ref{lemma: atomi}, $ At(\mathbf{R}_{\alpha(e)})=At(\mathbf{R})\setminus\{e\} $ and by induction hypothesis, $ \mathbf{R}_{\alpha(e)}=\prod_{c\in At(\mathbf{R}_{\alpha(e)})} \mathbf{R}_c $, whence our result follows.
\endproof 
\end{theorem}

\section{Appendix}\label{sec:app}

\noindent In section \ref{sec: 6}, we mentioned that the axiomatization of central element for the variety of integral involutive near semirings can be reduced to a minimal set of two identities. Indeed,  Theorem \ref{th: centrali nei near semirings} states that an element $ e $ of an involutive near semiring is central if and only if it satisfies the following identities: 
\begin{enumerate}
\item $ (e\cdot \alpha(x))+ (\alpha(e)\cdot \alpha(y))=\alpha((e\cdot x) + (\alpha(e)\cdot y)) $;
\item $ (e\cdot(x\cdot z)) + (\alpha(e)\cdot(y\cdot u))= ((e\cdot x)+ (\alpha(e)\cdot y))\cdot ((e\cdot z) + (\alpha(e)\cdot u)) $.
\end{enumerate}
\noindent
Here we provide a justification of the \emph{minimality} of this axiomatization. In fact, we will show in this section that identities (1) and (2) are independent.
\begin{example}\label{prop: near semiring che soddisfa (1) e non (2)}
\emph{The integral involutive near semiring $\mathbf{A} $, whose sum, multiplication and the antitone involution $ \alpha $ are defined in the following tables, satisfies (1) but not (2).}

\begin{table}[h]  \centering 
\begin{tabular}{r|r}
$\alpha $ & \\
 \hline 
 $0$ & $1$\\
 $ 1 $ & $ 0 $ \\
 $e$ & $a$ \\
 $ a $ & $ e $ \\
 $b$ & $c$ \\
 $c$ & $b$ 
\end{tabular} \hspace{.5cm}
\begin{tabular}{r|cccccc}
$+$ & $ 0 $ & $1$ & $e $ & $ a $ & $ b $ & $c$\\
\hline
    $ 0 $ & $ 0 $ & $1 $ & $e$ & $a$ & $a$ & $c$ \\
    $ 1 $ & $1$ & $1$ & $1$ & $1$ & $1$ & $1$ \\
    $ e $ & $ e $ & $1$ & $e$ & $1$ & $1$ & $e$ \\
    $ a $ & $a$ & $1$ & $1$ & $a$ & $a$ & $a$ \\
    $ b $ & $a$ & $1$ & $1$ & a & a & a \\
    $ c $ & $c$ & $1$ & $e$ & $a$ & $a$ & $c$
\end{tabular} \hspace{.5cm}
\begin{tabular}{r|rrrrrr}
$\cdot$ & $0$ & $1$ & $e$ & $a$ & $b$ & $c$\\
\hline
    $0$ & $0$ & $0$ & $0$ & $0$ & $0$ & $0$ \\
    $1$ & $0$ & $1$ & $e$ & $a$ & $b$ & $c$ \\
    $ e$ & $0$ & $e$ & $e$ & $0$ & $c$ & $c$ \\
    $a$ & $0$ & $a$ & $0$ & $a$ & $a$ & $0$ \\
    $b$ & $0$ & $b$ & $0$ & $a$ & $a$ & $0$ \\
    $c$ & $0$ & $c$ & $0$ & $0$ & $0$ & $0$
\end{tabular}
\end{table}
\end{example}
\noindent 
It is routine to check that $\mathbf{A} $ is an integral involutive near semiring, satisfying also identity (1). A counterexample to identity (2) is given by setting: $ x = b $, $ z = 1 $ and $ y = u = 0 $.
\begin{example}\label{contrex:due}
{\rm The integral involutive near semiring $\mathbf{B}$, whose sum, multiplication and the antitone involution $\alpha$ are defined in the following tables, satisfies (2) but not (1). \\

  \centering 
\begin{tabular}{r|r}
$\alpha$ & \\
\hline
$ 0$ & $1$\\
$ a$ & $a$\\ 
$1$ & $0$
\end{tabular} \hspace{.5cm}
\begin{tabular}{r|rrr}
$+$ & $0$ & $1$ & $a$\\
\hline
   $0$ & $0$ & $1$ & $a$ \\
   $1$ & $1$ & $1$ & $1$ \\
   $a$ & $a$ & $1$ & $a$
\end{tabular} \hspace{.5cm}
\begin{tabular}{r|rrr}
$\cdot$ & $0$ & $1$ & $a$\\
\hline
   $0$ & $0$ & $0$ & $a$ \\
   $1$ & $0$ & $1$ & $a$ \\
   $a$ & $0$ & $a$ & $a$
\end{tabular}\text{}}\\
\end{example}
\noindent
It is routine to check that $\mathbf B$ is an integral involutive near semiring satisfying equation (2). A counterexample to identity (1) is given by setting $e=0$ and $x=y=z=a$.

As a consequence of Examples \ref{prop: near semiring che soddisfa (1) e non (2)} and \ref{contrex:due} we conclude that

\begin{corollary}\label{cor:indpndt}
Equations $(1)$  and $(2)$ in Theorem \ref{th: centrali nei near semirings} are independent.
\end{corollary}

\section*{\textbf{Acknowledgement}}
\noindent{The first author acknowledges the Italian Ministry of Scientific Research (MIUR) for the support within the PRIN project ÒTheory of Rationality: logical, epistemological and computational aspectsÓ.
The work of the second author is supported by the bilateral Project ``New perspectives on residuated posets'' financed by Austrian Science Fund (FWF), project I 1923-N25, and the Czech Science Foundation (GA\v{C}R): project 15-34697L. The third author gratefully acknowledges the support of the Italian Ministry of Scientific
Research (MIUR) within the FIRB project ``Structures and Dynamics of Knowledge
and Cognition", Cagliari: F21J12000140001. Finally, we all thank Francesco Paoli for his valuable suggestions.
}


\end{document}